

\def\spaces{\space\space\space\space\space\space\space\space\space\space}
\def\spacess{\message{\spaces\spaces\spaces\spaces\spaces\spaces\spaces}}
\spacess
\spacess
\message{Annals of Mathematics Style: Current Version: 1.0. March 25, 1992}
\spacess
\spacess
\catcode`@=11



\font\elevensc=cmcsc10 scaled\magstephalf
\font\tensc=cmcsc10

\font\eightsc=cmcsc10 scaled800

\font\elevenrm=cmr10 scaled \magstephalf
\font\ninerm=cmr9
\font\eightrm=cmr8
\font\sixrm=cmr6
\font\fiverm=cmr5

\font\eleveni=cmmi10 scaled\magstephalf
\font\ninei=cmmi9
\font\eighti=cmmi8
\font\sixi=cmmi6
\font\fivei=cmmi5
\skewchar\ninei='177 \skewchar\eighti='177 \skewchar\sixi='177
\skewchar\eleveni= '177

\font\elevensy=cmsy10 scaled\magstephalf
\font\ninesy=cmsy9
\font\eightsy=cmsy8
\font\sixsy=cmsy6
\font\fivesy=cmsy5
\skewchar\ninesy='60 \skewchar\eightsy='60 \skewchar\sixsy='60
\skewchar\elevensy='60

\font\eighteenbf=cmbx10 scaled\magstep3

\font\twelvebf=cmbx10 scaled \magstep1
\font\elevenbf=cmbx10 scaled \magstephalf
\font\tenbf=cmbx10
\font\ninebf=cmbx9
\font\eightbf=cmbx8
\font\sixbf=cmbx6
\font\fivebf=cmbx5

\font\elevenit=cmti10 scaled\magstephalf
\font\nineit=cmti9
\font\eightit=cmti8

\font\eighteenmib=cmmi10 scaled \magstep3
\font\twelvemib=cmmib10 scaled \magstep1
\font\elevenmib=cmmib10 scaled\magstephalf
\font\tenmib=cmmib10
\font\eightmib=cmmib10 scaled 800 
\font\sixmib=cmmib10 scaled 600

\font\eighteensyb=cmsy10 scaled \magstep3
\font\twelvesyb=cmbsy10 scaled \magstep1
\font\elevensyb=cmbsy10 scaled \magstephalf
\font\tensyb=cmbsy10 
\font\eightsyb=cmbsy10 scaled 800
\font\sixsyb=cmbsy10 scaled 600
 
\font\elevenex=cmex10 scaled \magstephalf
\font\tenex=cmex10     
\font\eighteenex=cmex10 scaled \magstep3


\def\elevenpoint{\def\rm{\fam0\elevenrm}%
  \textfont0=\elevenrm \scriptfont0=\eightrm \scriptscriptfont0=\sixrm
  \textfont1=\eleveni \scriptfont1=\eighti \scriptscriptfont1=\sixi
  \textfont2=\elevensy \scriptfont2=\eightsy \scriptscriptfont2=\sixsy
  \textfont3=\elevenex \scriptfont3=\elevenex \scriptscriptfont3=\elevenex
  \def\bf{\fam\bffam\elevenbf}%
  \def\it{\fam\itfam\elevenit}%
  \textfont\bffam=\elevenbf \scriptfont\bffam=\eightbf
   \scriptscriptfont\bffam=\sixbf
\normalbaselineskip=13.9pt
  \setbox\strutbox=\hbox{\vrule height9.5pt depth4.4pt width0pt\relax}%
  \normalbaselines\rm}

\elevenpoint 

\def\ninepoint{\def\rm{\fam0\ninerm}%
  \textfont0=\ninerm \scriptfont0=\sixrm \scriptscriptfont0=\fiverm
  \textfont1=\ninei \scriptfont1=\sixi \scriptscriptfont1=\fivei
  \textfont2=\ninesy \scriptfont2=\sixsy \scriptscriptfont2=\fivesy
  \textfont3=\tenex \scriptfont3=\tenex \scriptscriptfont3=\tenex
  \def\it{\fam\itfam\nineit}%
  \textfont\itfam=\nineit
  \def\bf{\fam\bffam\ninebf}%
  \textfont\bffam=\ninebf \scriptfont\bffam=\sixbf
   \scriptscriptfont\bffam=\fivebf
\normalbaselineskip=11pt
  \setbox\strutbox=\hbox{\vrule height8pt depth3pt width0pt\relax}%
  \normalbaselines\rm}

\def\eightpoint{\def\rm{\fam0\eightrm}%
  \textfont0=\eightrm \scriptfont0=\sixrm \scriptscriptfont0=\fiverm
  \textfont1=\eighti \scriptfont1=\sixi \scriptscriptfont1=\fivei
  \textfont2=\eightsy \scriptfont2=\sixsy \scriptscriptfont2=\fivesy
  \textfont3=\tenex \scriptfont3=\tenex \scriptscriptfont3=\tenex
  \def\it{\fam\itfam\eightit}%
  \textfont\itfam=\eightit
  \def\bf{\fam\bffam\eightbf}%
  \textfont\bffam=\eightbf \scriptfont\bffam=\sixbf
   \scriptscriptfont\bffam=\fivebf
\normalbaselineskip=12pt
  \setbox\strutbox=\hbox{\vrule height8.5pt depth3.5pt width0pt\relax}%
  \normalbaselines\rm}


\def\eighteenbold{\def\rm{\fam0\eighteenrm}%
  \textfont0=\eighteenbf \scriptfont0=\twelvebf \scriptscriptfont0=\tenbf
  \textfont1=\eighteenmib \scriptfont1=\twelvemib\scriptscriptfont1=\tenmib
  \textfont2=\eighteensyb \scriptfont2=\twelvesyb\scriptscriptfont2=\tensyb
  \textfont3=\eighteenex \scriptfont3=\tenex \scriptscriptfont3=\tenex
  \def\bf{\fam\bffam\eighteenbf}%
  \textfont\bffam=\eighteenbf \scriptfont\bffam=\twelvebf
   \scriptscriptfont\bffam=\tenbf
\normalbaselineskip=20pt
  \setbox\strutbox=\hbox{\vrule height13.5pt depth6.5pt width0pt\relax}%
\everymath {\fam0 }
\everydisplay {\fam0 }
  \normalbaselines\bf}

\def\elevenbold{\def\rm{\fam0\elevenrm}%
  \textfont0=\elevenbf \scriptfont0=\eightbf \scriptscriptfont0=\sixbf%
  \textfont1=\elevenmib \scriptfont1=\eightmib \scriptscriptfont1=\sixmib%
  \textfont2=\elevensyb \scriptfont2=\eightsyb \scriptscriptfont2=\sixsyb%
  \textfont3=\elevenex \scriptfont3=\elevenex \scriptscriptfont3=\elevenex%
  \def\bf{\fam\bffam\elevenbf}%
  \textfont\bffam=\elevenbf \scriptfont\bffam=\eightbf%
   \scriptscriptfont\bffam=\sixbf%
\normalbaselineskip=14pt%
  \setbox\strutbox=\hbox{\vrule height10pt depth4pt width0pt\relax}%
\everymath {\fam0 }%
\everydisplay {\fam0 }%
  \normalbaselines\bf}


\hsize=31pc
\vsize=48pc

\parindent=22pt
\parskip=0pt

\widowpenalty=10000
\clubpenalty=10000

\topskip=12pt

\skip\footins=20pt
\dimen\footins=3in 

\abovedisplayskip=6.95pt plus3.5pt minus 3pt
\belowdisplayskip=\abovedisplayskip


\voffset=7pt
\def\bindingoffset{\ifodd\pageno \hoffset=-18pt\else\hoffset9pc\fi}

\newif\iftitle

\def\amheadline{\iftitle
\hbox to\hsize{\hss\currannalsline\hss}\else\line{\ifodd\pageno
\hfill\thetitle\hfill\llap{\elevenrm\folio}\else\rlap{\elevenrm\folio}
\hfill\theauthors\hfill\fi}\fi}

\headline={\amheadline}
\footline={\global\titlefalse}
\output={\bindingoffset\plainoutput}


\def\annalsline#1#2{\vfill\eject
\ifodd\pageno\else 
\line{\hfill}
\vfill\eject\fi
\global\titletrue
\def\currannalsline{\eightrm Annals of Mathematics,
{\eightbf#1} (#2), \thepages}}

\def\titleheadline#1{\def\one{#1}\ifx\one\empty\else
\gdef\thetitle{{\frenchspacing%
\let\\ \relax\eightsc\uppercase{#1}}}\fi}

\newif\ifshort

\let\shorttitle\titleheadline

\def\onpages#1#2{\def\thepages{#1--#2}}

\def\thismuchskip[#1]{\vskip#1pt}
\def\ilook{\ifx\next[ \let\go\thismuchskip\else
\let\go\relax\vskip1pt\fi\go}

\def\institution#1{\def\theinstitutions{\vbox{\baselineskip10pt
\def\and{{\eightrm and }}
\def\\{\futurelet\next\ilook}\eightsc #1}}}
\let\institutions\institution

\newwrite\auxfile

\def\startingpage#1{\def\one{#1}\ifx\one\empty\global\pageno=1\else
\global\pageno=#1\fi
\theoremcount=0 \eqcount=0 \sectioncount=0 
\openin1 \jobname.aux \ifeof1 
\onpages{#1}{???}
\else\closein1 \relax\input \jobname.aux
\onpages{#1}{\lastpage}
\fi\immediate\openout\auxfile=\jobname.aux
}

\def\endarticle{\ifRefsUsed\global\RefsUsedfalse%
\else\vskip21pt\theinstitutions%
\nobreak\vskip8pt\vbox{\thereceived\therevised}\fi%
\write\auxfile{\string\def\string\lastpage{\the\pageno}}}

\outer\def\bye{\endarticle\par \vfill \supereject \end}

\newif\ifacks
\long\def\acknowledgements#1{\def\one{#1}\ifx\one\empty\else
\global\ackstrue\vskip-\baselineskip%
\footnote{\ \unskip}{*#1}\fi}

\def\received#1{\def\thereceived{\medbreak
\centerline{\ninerm (Received #1)}}}

\def\revised#1{\def\therevised{\smallbreak
\centerline{\ninerm (Revised #1)}}}

\def\title#1{\vbox to76pt{\vfill
\baselineskip=18pt
\parindent=0pt
\overfullrule=0pt
\hyphenpenalty=10000
\everypar={\hskip\parfillskip\relax}
\hbadness=10000
\def\\ {\vskip1sp}
\eighteenbold#1\vskip1sp}
\titleheadline{#1}}

\newif\ifauthor
\def\author#1{\vskip15pt
\hbox to\hsize{\hss\tenrm By \tensc#1\ifacks\global\acksfalse*\fi\hss}
\ifshort\else\xdef\theauthors{{\eightsc\uppercase{#1}}}\fi%
\vskip7pt
\vskip\baselineskip
\global\authortrue\everypar={\global\authorfalse\everypar={}}}

\def\twoauthors#1#2{\vskip15pt
\hbox to\hsize{\hss%
\tenrm By \tensc#1 {\tenrm and} #2\ifacks\global\acksfalse*\fi\hss}
\ifshort\else\xdef\theauthors{{\eightsc\uppercase{#1 and #2}}}\fi
\vskip7pt
\vskip\baselineskip
\global\authortrue\everypar={\global\authorfalse\everypar={}}}


\newcount\theoremcount
\newcount\sectioncount
\newcount\eqcount

\newif\ifspecialnumon
\def\specialnumber#1{\global\specialnumontrue\gdef\spnum{#1}}

\def\eqnumber=#1 {\global\eqcount=#1 \global\advance\eqcount by-1\relax}
\def\sectionnumber=#1 {\global\sectioncount=#1 
\global\advance\sectioncount by-1\relax}
\def\proclaimnumber=#1 {\global\theoremcount=#1 
\global\advance\theoremcount by-1\relax}

\newif\ifsection
\newif\ifsubsection

\def\intro{\centerline{\bf Introduction}\global\everypar={}
\global\authorfalse\vskip6pt}

\def\elevenboldmath#1{$#1$\egroup}
\def\mathbold{\hbox\bgroup\elevenbold\elevenboldmath}

\def\section#1{\global\theoremcount=0
\global\eqcount=0
\ifauthor\global\authorfalse\else%
\vskip18pt plus 18pt minus 6pt\fi%
{\parindent=0pt\everypar={\hskip\parfillskip}%
\def\\ {\vskip1sp}\elevenpoint\bf%
\ifspecialnumon\global\specialnumonfalse$\rm\spnum$%
\gdef\sectnum{$\rm\spnum$}%
\else\interlinepenalty=10000%
\global\advance\sectioncount by1\relax\the\sectioncount%
\gdef\sectnum{\the\sectioncount}%
\fi.\hskip6pt#1\vrule width 0pt depth12pt}\hskip\parfillskip\break%
\global\sectiontrue%
\everypar={\global\sectionfalse\global\interlinepenalty=0\everypar={}}%
\ignorespaces}


\newif\ifspequation

\def\speqnu#1{\specialnumber{#1}\eqnu}

\let\eqno\leqno 

\newif\ifineqalignno
\let\saveleqalignno\leqalignno
\def\leqalignno{\let\eqnu\Eeqnu\saveleqalignno}

\let\eqalignno\leqalignno

\def\sectandeqnum{%
\ifspecialnumon\global\specialnumonfalse
$\rm\spnum$\gdef\eqnum{$\rm\spnum$}\else\global\firstlettertrue
\global\advance\eqcount by1 
\ifappend\applett\else\the\sectioncount\fi.%
\the\eqcount
\xdef\eqnum{\ifappend\applett\else\the\sectioncount\fi.\the\eqcount}\fi}

\def\eqnu{\leqno{\hbox{\elevenrm\ifspequation\else(\fi\sectandeqnum
\ifspequation\global\spequationfalse\else)\fi}}}

\def\Speqnu{\global\setbox\leqnobox=\hbox{\elevenrm
\ifspequation\else%
(\fi\sectandeqnum\ifspequation\global\spequationfalse\else)\fi}}

\def\Eeqnu{\hbox{\elevenrm
\ifspequation\else%
(\fi\sectandeqnum\ifspequation\global\spequationfalse\else)\fi}}

\newif\iffirstletter
\global\firstlettertrue
\def\eqletter#1{\global\specialnumontrue\iffirstletter\global\firstletterfalse
\global\advance\eqcount by1\fi
\gdef\spnum{\ifappend\applett\else\the\sectioncount\fi.\the\eqcount#1}\eqnu}

\newbox\leqnobox
\def\outsideeqnu#1{\global\setbox\leqnobox=\hbox{#1}}

\def\eatone#1{}

\def\dosplit#1#2{\vskip-.5\abovedisplayskip
\setbox0=\hbox{$\let\eqno\outsideeqnu%
\let\eqnu\Speqnu\let\leqno\outsideeqnu#2$}%
\setbox1\vbox{\noindent\hskip\wd\leqnobox\ifdim\wd\leqnobox>0pt\hskip1em\fi%
$\displaystyle#1\mathstrut$\hskip0pt plus1fill\relax
\vskip1pt
\line{\hfill$\let\eqnu\eatone\let\leqno\eatone%
\displaystyle#2\mathstrut$\ifmathqed~~\qed\fi}}%
\copy1
\ifvoid\leqnobox
\else\dimen0=\ht1 \advance\dimen0 by\dp1
\vskip-\dimen0
\vbox to\dimen0{\vfill
\hbox{\unhbox\leqnobox}
\vfill}
\fi}



\def\testforbreak#1\splitmath #2@{\def\mathtwo{#2}\ifx\mathtwo\empty%
#1$$%
\ifmathqed\vskip-\belowdisplayskip
\setbox0=\vbox{\let\eqno\relax\let\eqnu\relax$\displaystyle#1$}%
\vskip-\ht0\vskip-3.5pt\hbox to\hsize{\hfill\qed}
\vskip\ht0\vskip3.5pt\fi
\else$$\vskip-\belowdisplayskip
\vbox{\dosplit{#1}{\let\eqno\eatone
\let\splitmath\relax#2}}%
\nobreak\vskip.5\belowdisplayskip
\noindent\ignorespaces\fi}


\newif\ifmathqed



\newcount\linenum
\newcount\colnum
\def\spline{\omit&\multispan{\the\colnum}{\hrulefill}\cr}
\def\colcounter{\ifnum\linenum=1\global\advance\colnum by1\fi}

\def\everyline{\noalign{\global\advance\linenum by1\relax}%
\ifnum\linenum=2\spline\fi}

\def\mtable{\bgroup\offinterlineskip
\everycr={\everyline}\global\linenum=0
\halign\bgroup\vrule height 10pt depth 4pt width0pt
\hfill$##$\hfill\hskip6pt\ifnum\linenum>1
\vrule\fi&&\colcounter\hskip12pt\hfill$##$\hfill\hskip12pt\cr}

\def\endmtable{\crcr\egroup\egroup}




\def\xast{*}
\newcount\intable
\newcount\mathcol
\newcount\savemathcol
\newcount\topmathcol
\newdimen\arrayhspace
\newdimen\arrayvspace

\arrayhspace=8pt 
\arrayvspace=12pt 

\newif\ifdollaron

\def\mathalign#1{\def\arg{#1}\ifx\arg\xast%
\let\go\relax\else\let\go\mathalign%
\global\advance\mathcol by1 %
\global\advance\topmathcol by1 %
\expandafter\def\csname  mathcol\the\mathcol\endcsname{#1}%
\fi\go}

\def\arraypickapart#1]#2*{\if#1c \ifmmode\vcenter\else
\global\dollarontrue$\vcenter\fi\else%
\if#1t\vtop\else\if#1b\vbox\fi\fi\fi\bgroup%
\def\one{#2}}

\def\arraystrut{\vrule height .7\arrayvspace depth .3\arrayvspace width 0pt}

\def\array#1#2*{\def\firstarg{#1}%
\if\firstarg[ \def\two{#2} \expandafter\arraypickapart\two*\else%
\ifmmode\vcenter\else\vbox\fi\bgroup \def\one{#1#2}\fi%
\global\everycr={\noalign{\global\mathcol=\savemathcol\relax}}%
\def\\ {\cr}%
\global\advance\intable by1 %
\ifnum\intable=1 \global\mathcol=0 \savemathcol=0 %
\else \global\advance\mathcol by1 \savemathcol=\mathcol\fi%
\expandafter\mathalign\one*%
\mathcol=\savemathcol %
\halign\bgroup&\hskip.5\arrayhspace\arraystrut%
\global\advance\mathcol by1 \relax%
\expandafter\if\csname mathcol\the\mathcol\endcsname r\hfill\else%
\expandafter\if\csname mathcol\the\mathcol\endcsname c\hfill\fi\fi%
$\displaystyle##$%
\expandafter\if\csname mathcol\the\mathcol\endcsname r\else\hfill\fi\relax%
\hskip.5\arrayhspace\cr}

\def\endarray{\crcr\egroup\egroup%
\global\mathcol=\savemathcol %
\global\advance\intable by -1\relax%
\ifnum\intable=0 %
\ifdollaron\global\dollaronfalse $\fi
\loop\ifnum\topmathcol>0 %
\expandafter\def\csname  mathcol\the\topmathcol\endcsname{}%
\global\advance\topmathcol by-1 \repeat%
\global\everycr={}\fi%
}

\def\big#1{{\hbox{$\left#1\vbox to 10pt{}\right.\n@space$}}}
\def\Big#1{{\hbox{$\left#1\vbox to 13pt{}\right.\n@space$}}}
\def\bigg#1{{\hbox{$\left#1\vbox to 16pt{}\right.\n@space$}}}
\def\Bigg#1{{\hbox{$\left#1\vbox to 19pt{}\right.\n@space$}}}


\def\figcaption#1#2#3{\topinsert
\vskip4pt 
\vbox to#3{\vfill}\vskip1sp
\setbox0=\hbox{\eightsc Figure #1.\hskip12pt\eightpoint\rm #2}
\ifdim\wd0>\hsize
\noindent\eightsc Figure #1.\hskip12pt\eightpoint\rm #2\else
\centerline{\eightsc Figure #1.\hskip12pt\eightpoint\rm #2}\fi
\vskip16pt
\endinsert}

\def\wfig#1#2#3{\topinsert
\vskip4pt 
\hbox to\hsize{\hss\vbox{\hrule height .25pt width #3
\hbox to #3{\vrule width .25pt height #2\hfill\vrule width .25pt height #2}
\hrule height.25pt}\hss}
\vskip1sp
\centerline{\eightsc Figure #1}
\vskip16pt
\endinsert}

\def\wfigcaption#1#2#3#4{\topinsert
\vskip4pt 
\hbox to\hsize{\hss\vbox{\hrule height .25pt width #4
\hbox to #4{\vrule width .25pt height #3\hfill\vrule width .25pt height #3}
\hrule height.25pt}\hss}
\vskip1sp
\setbox0=\hbox{\eightsc Figure #1.\hskip12pt\eightpoint\rm #2}
\ifdim\wd0>\hsize
\noindent\eightsc Figure #1.\hskip12pt\eightpoint\rm #2\else
\centerline{\eightsc Figure #1.\hskip12pt\eightpoint\rm #2}\fi
\vskip16pt
\endinsert}

\def\tabcaption#1#2{\vskip6pt
\setbox0=\hbox{\eightsc Table #1.\hskip12pt #2}
\ifdim\wd0>\hsize
\noindent\eightsc Table #1.\hskip12pt\eightpoint #2
\else
\centerline{\eightsc Table #1.\hskip12pt\eightpoint #2}
\fi
\vskip6pt}

\def\endinsert{\egroup\if@mid\dimen@\ht\z@\advance\dimen@\dp\z@ 
\advance\dimen@ 12\p@\advance\dimen@\pagetotal\ifdim\dimen@ >\pagegoal 
\@midfalse\p@gefalse\fi\fi\if@mid\smallskip\box\z@\bigbreak\else
\insert\topins{\penalty 100 \splittopskip\z@skip\splitmaxdepth\maxdimen
\floatingpenalty\z@\ifp@ge\dimen@\dp\z@\vbox to\vsize {\unvbox \z@ 
\kern -\dimen@ }\else\box\z@\nobreak\smallskip\fi}\fi\endgroup}

\def\pagecontents{
\ifvoid\topins \else\iftitle\else 
\unvbox \topins \fi\fi \dimen@ =\dp \@cclv \unvbox 
\@cclv 
\ifvoid\topins\else\iftitle\unvbox\topins\fi\fi
\ifvoid \footins \else \vskip \skip \footins \footnoterule 
\unvbox \footins \fi \ifr@ggedbottom \kern -\dimen@ \vfil \fi}


\newif\ifappend

\def\appendix#1#2{\def\applett{#1}\def\two{#2}%
\global\appendtrue
\global\theoremcount=0
\global\eqcount=0
\vskip18pt plus 18pt
\vbox{\parindent=0pt
\everypar={\hskip\parfillskip}
\def\\ {\vskip1sp}\elevenpoint\bf Appendix%
\ifx\applett\empty\gdef\applett{A}\ifx\two\empty\else.\fi%
\else\ #1.\fi\hskip6pt#2\vskip12pt}%
\global\sectiontrue%
\everypar={\global\sectionfalse\everypar={}}\nobreak\ignorespaces}

\newif\ifRefsUsed
\long\def\references{\global\RefsUsedtrue\vskip21pt
\theinstitutions
\everypar={}\global\bibnum=0
\vskip20pt\goodbreak\bgroup
\vbox{\centerline{\eightsc References}\vskip6pt}%
\ifdim\maxbibwidth>0pt\leftskip=\maxbibwidth\else\leftskip=18pt\fi%
\ifdim\maxbibwidth>0pt\parindent=-\maxbibwidth\else\parindent=-18pt\fi%
\everypar={\amref}%
\ninepoint
\frenchspacing
\nobreak\ignorespaces}

\def\endreferences{\vskip1sp\egroup\everypar={}%
\nobreak\vskip8pt\vbox{\thereceived\therevised}}

\newcount\bibnum

\def\amref#1 {\def\one{#1}\global\advance\bibnum by1%
\immediate\write\auxfile{\string\expandafter\string\def\string\csname
\space #1\string\endcsname{[\the\bibnum]}}%
\leavevmode\hbox to18pt{\hbox to13.2pt{\hss[\the\bibnum]}\hfill}}

\def\bibline{\hbox to30pt{\hrulefill}\/\/}

\def\name#1{{\eightsc#1}}

\newdimen\maxbibwidth
\def\AuthorRefNames [#1] {\def\amref{\spamref}
\setbox0=\hbox{[#1] }\global\maxbibwidth=\wd0\relax}

\def\spamref[#1] {\leavevmode\hbox to\maxbibwidth{\hss[#1]\hfill}}


\def\footnoterule{\kern-3pt\hrule width1in height.5pt\kern2.5pt}

\let\savefootnote\footnote
\def\footnote#1#2{%
\savefootnote{#1}{{\eightpoint\normalbaselineskip11pt
\normalbaselines#2}}}

\def\vfootnote#1{%
\insert \footins \bgroup \eightpoint\normalbaselineskip11pt\normalbaselines
\interlinepenalty \interfootnotelinepenalty
\splittopskip \ht \strutbox \splitmaxdepth \dp \strutbox \floatingpenalty 
\@MM \leftskip \z@skip \rightskip \z@skip \spaceskip \z@skip 
\xspaceskip \z@skip
{#1}$\,$\footstrut \futurelet \next \fo@t}


\newif\iffirstadded
\newif\ifadded

\def\addedlett{}

\def\alltheoremnums{%
\ifspecialnumon\global\specialnumonfalse
\ifadded\global\addedfalse
\iffirstadded\global\firstaddedfalse
\global\advance\theoremcount by1 \fi
\ifappend\applett\else\the\sectioncount\fi.\the\theoremcount\addedlett%
\xdef\theoremnum{\ifappend\applett\else\the\sectioncount\fi.%
\the\theoremcount\addedlett}%
\else$\rm\spnum$\def\theoremnum{$\rm\spnum$}\fi%
\else\global\firstaddedtrue
\global\advance\theoremcount by1 
\ifappend\applett\else\the\sectioncount\fi.\the\theoremcount%
\xdef\theoremnum{\ifappend\applett\else\the\sectioncount\fi.%
\the\theoremcount}\fi}

\def\allcorolnums{%
\ifspecialnumon\global\specialnumonfalse
\ifadded\global\addedfalse
\iffirstadded\global\firstaddedfalse
\global\advance\corolcount by1 \fi
\the\corolcount\addedlett%
\else$\rm\spnum$\def\corolnum{$\rm\spnum$}\fi%
\else\global\advance\corolcount by1 
\the\corolcount\fi}


\newcount\corolcount
\def\xcorol{Corollary}
\def\xtheorem{Theorem}
\def\xmaintheorem{Main Theorem}

\newif\ifthtitle

\let\saverparen)
\let\savelparen(
\def\rmparenl{{\rm(}}
\def\rmparenr{{\rm\/)}}
{
\catcode`(=13
\catcode`)=13
\gdef\makeparensRM{
\catcode`(=13
\catcode`)=13
\let(=\rmparenl
\let)=\rmparenr
\everymath{\let(\savelparen
\let)\saverparen}
\everydisplay{\let(\savelparen
\let)\saverparen}}
}

\medskipamount=8pt plus.1\baselineskip minus.05\baselineskip
\def\proclaim#1{\vskip-\lastskip
\def\one{#1}\ifx\one\xtheorem\global\corolcount=0\fi
\ifsection\global\sectionfalse\vskip-6pt\fi
\medskip
{\elevensc#1}%
\ifx\one\xmaintheorem\global\corolcount=0
\gdef\theoremnum{Main Theorem}\else%
\ifx\one\xcorol\ \allcorolnums\else\ \alltheoremnums\fi\fi%
\ifthtitle\ \global\thtitlefalse{\rm(\thethtitle)}\fi.%
\hskip.5pc\bgroup\makeparensRM
\it\ignorespaces}

\def\nonumproclaim#1{\vskip-\lastskip
\def\one{#1}\ifx\one\xtheorem\global\corolcount=0\fi
\ifsection\global\sectionfalse\vskip-6pt\fi
\medskip
{\elevensc#1}.\ifx\one\xmaintheorem\global\corolcount=0
\gdef\theoremnum{Main Theorem}\fi\hskip.5pc%
\bgroup\makeparensRM\it\ignorespaces}

\def\endproclaim{\egroup\medskip}


\def\xproof{Proof}
\def\xremark{Remark}
\def\xcase{Case}
\def\xsubcase{Subcase}
\def\xconjecture{Conjecture}
\def\xstep{Step}
\def\xof{of}

\def\deconstruct#1 #2 #3 #4 #5 @{\def\one{#1}\def\two{#2}\def\three{#3}%
\def\four{#4}%
\ifx\two\empty #1\else%
\ifx\one\xproof%
\ifx\two\xof%
  \ifx\three\xcorol Proof of Corollary \rm#4\else%
     \ifx\three\xtheorem Proof of Theorem \rm#4\else\xone\fi%
  \fi\fi%
\else\xone\fi\fi.}

\def\pickup#1 {\def\this{#1}%
\ifx\this\xproof\global\let\go\demoproof
\global\let\enddemo\endproof\else
\ifx\this\xremark\global\let\go\demoremark\else
\ifx\this\xcase\global\let\go\demostep\else
\ifx\this\xsubcase\global\let\go\demostep\else
\ifx\this\xconjecture\global\let\go\demostep\else
\ifx\this\xstep\global\let\go\demostep\else
\global\let\go\demoproof\fi\fi\fi\fi\fi\fi}

\def\demo#1{\vskip-\lastskip
\ifsection\global\sectionfalse\vskip-6pt\fi
\def\one{#1 }\def\two{#1*}%
\setbox0=\hbox{\expandafter\pickup\one}\expandafter\go\two}

\def\numbereddemo#1{\vskip-\lastskip
\ifsection\global\sectionfalse\vskip-6pt\fi
\def\two{#1*}%
\expandafter\demoremark\two}

\def\demoproof#1*{\medskip\def\xone{#1}
{\ignorespaces\it\expandafter\deconstruct\xone {} {} {} {} {} @%
\unskip\hskip6pt}\rm\ignorespaces}

\def\demoremark#1*{\medskip
{\it\ignorespaces#1\/} \alltheoremnums\unskip.\hskip1pc\rm\ignorespaces}

\def\demostep#1 #2*{\vskip4pt
{\it\ignorespaces#1\/} #2.\hskip1pc\rm\ignorespaces}

\def\enddemo{\medskip}

\def\endproof{\ifmathqed\global\mathqedfalse\medskip\else
\parfillskip=0pt~~\hfill\qed\medskip
\fi\global\parfillskip0pt plus 1fil\relax
\gdef\enddemo{\medskip}}

\def\qed{\vbox{\hrule\hbox{\vrule height6pt\hskip6pt\vrule}\hrule}}








\def\stripbs#1#2*{\def\one{#2}}
\def\newdef#1#2#{\expandafter\stripbs\string#1*
\defining{#1}{#2}}

\def\defining#1#2#3{\ifdefined{\one}{\message{<<<Sorry, 
\string#1\space
is already defined. Please use a new name.>>>}}
{\expandafter\gdef#1#2{#3}}}

\def\emptyspace{ }
\def\nextthing{}
\def\newline{***}
\def\eatone#1{ }

\def\lookatspace#1{\ifcat\noexpand#1\ \else%
\gdef\nextthing{}\xdef\next{#1}%
\ifx\next\emptyspace%
\let\nextthing\emptyspace\else\ifx\next\newline%
\gdef\nextthing{\eatone}\fi\fi\fi\egroup\nextthing#1}

{\catcode`\^^M=\active%
\gdef\spacer{\bgroup\catcode`\^^M=\active%
\let^^M=\newline\obeyspaces\lookatspace}}

\def\ref#1{\seeifdefined{#1}#1\spacer}


\def\seeifdefined#1{\expandafter\stripbs\string#1*%
\crorefdefining{#1}}

\newif\ifcromessage
\global\cromessagetrue

\def\crorefdefining#1{\ifdefined{\one}{}
{\ifcromessage\global\cromessagefalse%
\message{\spaces\spaces\spaces\spaces\spaces\spaces\spaces}%
\message{<Undefined reference.}%
\message{Please TeX file once more to have accurate cross-references.>}%
\message{\spaces\spaces\spaces\spaces\spaces\spaces\spaces}\fi%
\expandafter\gdef#1{??}}}

\def\label#1#2*{\gdef\ctest{#2}\seeifdefined{#1}%
\ifx\empty\ctest\write\auxfile{\noexpand\def\noexpand#1{\tocstuff}}%
\else\def\ctemp{#2}%
\immediate\write\auxfile{\noexpand\def\noexpand#1{\ctemp}}\fi}

\def\ifdefined#1#2#3{\expandafter\ifx\csname#1\endcsname\relax#3\else#2\fi}




\def\articlecontents{
\vskip20pt\centerline{\bf Table of Contents}\everypar={}\vskip6pt
\bgroup \leftskip=3pc \parindent=-2pc 
\def\item##1{\vskip1sp\indent\hbox to2pc{##1.\hfill}}}

\def\endcontents{\vskip1sp\leftskip=0pt\egroup}

\def\journalcontents{\vfill\eject
\def\currannalsline{\hfill}
\global\titletrue
\vglue3.5pc
\centerline{\tensc\hskip12pt TABLE OF CONTENTS}\everypar={}\vskip30pt
\bgroup \leftskip=34pt \rightskip=-12pt \parindent=-22pt 
  \def\\ {\vskip1sp\noindent}
\def\pagenum##1{\unskip\parfillskip=0pt\dotfill##1\vskip1sp
\parfillskip=0pt plus 1fil\relax}
\def\name##1{{\tensc##1}}}

\def\indexcontents{\vfill\eject
\leftline{\tensc \quad Annals of Mathematics}
\vglue-2truept
\leftline{\quad\phantom{AN}\ninerm   Vol.\ 148, Nos.\ 1, 2, 3}
\vglue-2truept
\leftline{\quad\phantom{ANNALS} \ninerm 1998}
\def\currannalsline{\hfill} 
\def\amheadline{\iftitle
\hbox to\hsize{\hss\currannalsline\hss}\else\line{\ifodd\pageno
\hfill\thetitle\hfill\llap{\elevenrm\folio}\else\rlap{\elevenrm\folio}
\hfill\fi}\fi}

\headline={\amheadline}
\global\titletrue
\vglue2.5pc 
\centerline{\tensc\hskip12pt INDEX}\everypar={}\vskip18pt
\bgroup \leftskip=34pt \rightskip=-12pt \parindent=-22pt 
  \def\\ {\vskip1sp\noindent}
\def\pagenum##1{\unskip\parfillskip=0pt\dotfill##1\vskip1sp
\parfillskip=0pt plus 1fil\relax}
\def\name##1{{\tensc##1}}}


\institution{}
\onpages{0}{0}
\def\lastpage{???}
\def\thetitle{Title ???}
\def\theauthors{Authors ???}
\def\thereceived{}
\def\therevised{}

\catcode`\@=12

\annalsline{158}{2003}
\received{January 30, 2001}
\revised{July 22, 2002}
 \startingpage{473}
 
\catcode`\@=11
\font\twelvemsb=msbm10 scaled 1100

\font\ninemsb=msbm10 scaled 800
\newfam\msbfam
\textfont\msbfam=\twelvemsb  \scriptfont\msbfam=\ninemsb
  \scriptscriptfont\msbfam=\ninemsb
\def\msb@{\hexnumber@\msbfam}
\def\Bbb{\relax\ifmmode\let\next\Bbb@\else
 \def\next{\errmessage{Use \string\Bbb\space only in math
mode}}\fi\next}
\def\Bbb@#1{{\Bbb@@{#1}}}
\def\Bbb@@#1{\fam\msbfam#1}
\catcode`\@=12

 \font\tenrm=cmr10
\def\ritem#1{\vglue4pt \item{\rm #1}}
\def\rritem#1{\vglue4pt\itemitem{\rm #1}}
 \def\R{{\Bbb R}}
\def\N{{\Bbb N}}
\def\C{{\Bbb C}}
\def\P{{\Bbb P}}
\def\Q{{\Bbb Q}}
\def\Z{{\Bbb Z}}
\def\L{{\Bbb L}}
\def\K{{\Bbb K}}
\def\T{{\Bbb T}}
\phantom{what}

\phantom{why}

\title{Hypersurface complements, Milnor fibers\\ and higher homotopy groups\\ of arrangments} 
\shorttitle{Hypersurface complements} 
 \twoauthors{Alexandru Dimca}{Stefan Papadima}
 \institutions{Laboratoire d'Analyse et Geometrie, UMR 5467 CNRS, Universit\'e Bordeaux I, 33405 Talence, France
 \\
{\eightpoint {\it E-mail address\/}: 
 dimca@math.u-bordeaux.fr}
\\ \vglue6pt
Institute of Mathematics of the Romanian Academy, Bucharest, Romania\\
{\eightpoint {\it E-mail address\/}:
 spapadim@stoilow.imar.ro}
}
 
\vfil
\intro

\vfil   

The interplay between geometry and topology on complex algebraic
varieties
is a classical theme that goes back to Lefschetz [L] and Zariski [Z]
and is always present on the scene; see for instance the work by Libgober
[Li].
In this paper we study complements of
hypersurfaces, with   special attention to the case of hyperplane
arrangements as discussed in Orlik-Terao's book [OT1].

Theorem 1 expresses the degree of the gradient map associated to any
homogeneous polynomial $h$ as the number of $n$-cells that have to be
added to a generic hyperplane section $D(h) \cap H$ to obtain the
complement in $\P^n$, $D(h)$, of the projective hypersurface $V(h)$.
Alternatively, by
results of L\^e [Le2] one knows that the affine piece $V(h)_a=V(h)
\setminus H$ of $V(h)$ has the homotopy type of a bouquet of $(n-1)$-spheres.
Theorem 1 can then be restated by saying that the degree of the
gradient map coincides with the number of these $(n-1)$-spheres.
In this form, our result is reminiscent of Milnor's equality between
the  degree of the local gradient map and the number of spheres in the
Milnor fiber associated to an  isolated hypersurface singularity [M].

This topological description of the degree of the gradient map has as a
direct consequence a positive answer to a conjecture by Dolgachev [Do]
on polar Cremona transformations; see Corollary 2. Corollary 4 and the
end of Section 3 contain stronger versions of some of the results in
[Do] and some related matters.

Corollary 6 (obtained
independently by  Randell [R2,3]) reveals a striking feature of complements of
hyperplane arrangements. They possess
a {\it minimal} CW-structure, i. e., a CW-decomposition with
exactly as many $k$-cells as the $k$-th Betti number, for all $k$.
Minimality may be viewed as an improvement of the Morse inequalities
for twisted homology (the main result of Daniel Cohen in [C]), from
homology to the level of cells; see Remark 12  (ii).\eject

In the second part of
our paper, we investigate the higher homotopy groups
of complements of complex hyperplane arrangements
(as $\pi_1$-modules).
By the classical work of Brieskorn [B] and Deligne [De], it is known
that such a complement is often aspherical. The first explicit computation
of nontrivial homotopy groups of this type has been performed by Hattori [Hat],
in 1975. This remained the only example of this kind, until [PS] was published.

Hattori proved that, up to homotopy, the complement of a general position arrangement
is a skeleton of the standard minimal CW-structure of a torus.
From this, he derived a free resolution of the first nontrivial higher homotopy group.
We use the techniques developed in the first part of our paper
to generalize Hattori's homotopy type formula, for all sufficiently generic sections
of aspherical arrangements (a framework inspired from the stratified
Morse theory of Goresky-MacPherson [GM]); see Proposition 14.
Using the approach by minimality from [PS], we can 
to generalize the Hattori presentation in Theorem 16, and the
Hattori resolution in Theorem 18. The above framework provides a unified treatment
of all explicit computations related to nonzero higher
homotopy groups of arrangements available in the literature, to the best of our knowledge.
It also gives examples exhibiting a nontrivial homotopy group,
$\pi_q$, for {\it all} $q$; see the end of Section 5.

The associated {\it combinatorics} plays an important role in
arrangement theory. By `combinatorics' we mean the pattern of
intersection of the hyperplanes, encoded by the associated
intersection lattice. For instance, one knows, by the work of Orlik-Solomon [OS],
that the cohomology ring of the complement is determined by
the combinatorics. On the other hand, the examples of Rybnikov [Ry]
show that $\pi_1$ is not combinatorially determined, in general.
One of the most basic questions in the field is to identify
the precise amount of topological information on the complement
that is determined by the combinatorics.

In Corollary 21 we consider the associated graded chain complex,
with respect to the $I$-adic filtration of the group ring $\Z \pi_1$,
of the $\pi_1$-equivariant chain complex of the universal cover,
constructed from an arbitrary minimal CW-structure of any
arrangement complement. We prove that the associated graded is
always combinatorially determined, with $\Q$-coefficients,
and that this actually holds over $\Z$, for the class of
{\it hypersolvable} arrangements introduced in [JP1].
We deduce these properties from a general result, namely
Theorem~20, where we show that the associated graded
equivariant chain complex of the universal cover of a
minimal CW-complex, whose cohomology ring is generated in degree one,
is determined by $\pi_1$ and the cohomology ring.
\pagegoal=50pc

There is a rich supply of examples which fit into our framework of
generic sections of aspherical arrangements. Among them,
we present in Theorem 23 a large class of combinatorially defined
hypersolvable examples, for which the associated graded module
of the first higher nontrivial homotopy group of the complement
is also combinatorially determined.
\eject

\section{The main results}

There is a gradient map associated to any nonconstant homogeneous
polynomial
$h \in \C[x_0,\ldots ,x_n]$ of degree $d$, namely
$$  {\rm grad}(h): D(h) \to \P^n, ~~~~(x_0:\cdots :x_n) \mapsto
(h_0(x):\cdots :h_n(x))$$
where $D(h)= \{x \in \P^n | h(x) \not= 0\}$ is the principal open set
associated to $h$ and $h_i={\partial h \over \partial x_i}$. This map
corresponds to the polar Cremona transformations considered by Dolgachev
in [Do].
Our first  result is the following topological description of the
degree
of the gradient map $  {\rm grad}(h)$.

\specialnumber{1}\proclaim{Theorem}   For any nonconstant homogeneous
polynomial
$h \in\break \C[x_0,\ldots ,x_n]$,
 the complement $D(h)$ is homotopy equivalent to a {\rm CW} complex obtained
 from $D(h) \cap H$ by attaching ${\rm deg}({\rm grad}(h)) $ cells of dimension
$n${\rm ,}
 where $H$ is a generic hyperplane
in $\P^n$. In particular{\rm ,} one has
$${\rm deg}({\rm grad}(h))= (-1)^n \chi (D(h) \setminus H).$$
\endproclaim

Note that the meaning of `generic' here is quite explicit:
the hyperplane $H$ has to be transversal to a stratification
of the projective hypersurface $V(h)$ defined by $h=0$ in $\P^n$.

The Euler characteristic in the above statement can be replaced by a
Betti number as follows. As noted in the introduction,
the affine part $V(h)_a=V(h) \setminus H$ of $V(h)$ has the homotopy
type of a bouquet of $(n-1)$-spheres. Using the additivity
of the Euler characteristic with respect to constructible partitions we
get
$${\rm deg}({\rm grad}(h))=b_{n-1}(V(h)_a).$$
 We use this form of Theorem 1 at the end of Section 3 to give a
topological easy proof of Theorem 4 in [Do].

 Theorem 1 looks similar to a conjecture on hyperplane arrangements by
Varchenko [V] proved by Orlik and Terao [OT3] and by Damon [Da2], but
in fact it is not; see our discussion after Lemma 5.

On the other hand, with  some transversality conditions for the
irreducible factors of $h$, Damon has obtained a local form of Theorem 1 in
which
$$D(h) \setminus H = \{ x \in \C^{n+1}| \ell (x)=1 \} \setminus \{ x
\in \C^{n+1}| h (x)=0 \}$$
($\ell (x)=0$ being an equation for $H$) is replaced by
$$V_t  \setminus \{ x \in \C^{n+1}| h (x)=0 \}$$
with $V_t$ the Milnor fiber of an isolated complete intersection
singularity $V$ at the origin of $\C^{n+1}$; see
[Da2, Th.~1]. In such a situation the corresponding Euler number is
explicitly computed as a sum of the Milnor number $\mu (V)$ and a
``singular Milnor number''; see [Da2, Th.~1] and [Da3, Th.~1 or Cor.~2].

\nonumproclaim{{C}orollary 2}  The degree of the gradient map
${\rm grad}(h)$ depends only on the reduced polynomial $h_r$ associated to $h$.
\endproclaim

This gives a positive answer to Dolgachev's conjecture at the end of
Section~3
in [Do], and it follows directly from Theorem 1, since $D(h)=D(h_r)$.
 
Let $f \in \C[x_0,\ldots ,x_n]$ be a homogeneous polynomial of degree $e>0$
with
global Milnor fiber $F=\{ x \in \C^{n+1}| f(x)=1\}$; see for instance
[D1] for more on such varieties. Let $g: F \setminus N \to \R$
be the function $g(x)=h(x) {\overline h}(x)$,
where $N=\{ x \in \C^{n+1}| h(x)=0\}$. Then we have the following:

\nonumproclaim{Theorem 3} \hskip-8pt  For any reduced homogeneous polynomial
$h \in \C[x_0,\ldots ,x_n]$ there is a Zariski open and dense subset ${\cal
U}$
in the space of homogeneous polynomials of degree $e>0$ such that for
any $f \in {\cal U}$
one has the following\/{\rm :}\/

\ritem{(i)} the function $g$ is a Morse function{\rm ;}

\ritem{(ii)} the Milnor fiber $F$ is homotopy equivalent to a CW complex
obtained from $F \cap N$ by attaching $|C(g)|$ cells
of dimension $n${\rm ,} where $C(g)$ is the critical set of the Morse
function $g${\rm ;}

\ritem{(iii)} the intersection $F \cap N$ is homotopy equivalent to a bouquet
of $|C(g)|-(e-1)^{n+1}$ spheres $S^{n-1}$. 

\endproclaim
 
In some cases the open set ${\cal U}$ can be explicitly described, as 
in Corollary~7 below. In general this task is a difficult one in view
of the proof of Theorem~3.
The  claim (iii) above, in the special case $e=1$, gives a new proof
for L\^e's result mentioned in the introduction.

\vglue6pt {\elevensc Lefschetz Theorem on generic hyperplane complements in\break hypersurfaces}.
{\it For any projective hypersurface $V(h): h=0$ in $\P ^n$ and any generic
hyperplane $H$ in $\P ^n$ the affine hypersurface given by the
complement $V(h) \setminus H$ is homotopy equivalent to a bouquet of  spheres
$S^{n-1}$.}

\vglue6pt

We point out that both Theorem 1 and Theorem 3 follow from the
results by Hamm in [H]. In the case of Theorem 1, the homotopy-type
claim is a direct consequence of [H, Th.~5], the new part being the
relation between the number of $n$-cells and the degree of the
gradient map ${\rm grad}(h)$. We establish this equality by
using polar curves and complex Morse theory; see Section 2.

On the other hand, in Theorem 3 the main claim is that concerning
the homotopy-type and this follows from  [H, Prop.~3], by a geometric
argument described in Section 3 and involving a key result by Hironaka.
An alternative proof may also be given using Damon's work [Da1, Prop.~9.14], a result which extends
previous results by Siersma [Si] and Looijenga [Lo].

Our results above have interesting implications for the topology of
hyperplane arrangements which were our initial
motivation in this study.
Let ${\cal A}$ be a hyperplane arrangement in the complex projective
space $\P^n$, with $n>0$. Let $d>0$ be the number of hyperplanes in
this arrangement and choose a linear equation $ \ell_i(x)=0$ for
each hyperplane $H_i$ in ${\cal A}$,
for $i=1,\ldots ,d.$

Consider the homogeneous polynomial
$Q(x)= \prod _{i=1}^d \ell_i(x) \in \C[x_0,\ldots ,x_n]$ and the
corresponding
principal open set $M=M({\cal A})=D(Q)=\P^n \setminus \mathbold{\cup}_{i=1}^d
H_i$.
The topology of the hyperplane arrangement complement $M$ is a central
object of study in the theory of hyperplane arrangements, see
Orlik-Terao [OT1].
As a  consequence of Theorem 1 we prove the following:

\nonumproclaim{{C}orollary 4} {\rm (1)}   For any projective arrangement
${\cal A}$
as above one has
$$b_n(D(Q))={\rm deg}({\rm grad}(Q)).$$ 

{\rm (2)}  In particular\/{\rm :}\/

\ritem{(a)} The following are equivalent\/{\rm :}\/

\rritem{(i)} the morphism ${\rm grad}(Q)$ is dominant{\rm ;}

\rritem{(ii)}  $b_n(D(Q)) >0${\rm ;}

\rritem{(iii)} the projective arrangement ${\cal A}$ is  essential{\rm ;}
 i.e.{\rm ,} the intersection $\mathbold{\cap}_{i=1}^d H_i\!$ is empty.

\ritem{(b)} If $b_n(D(Q))>0$ then $d \leq n+b_n(D(Q))$. As special cases\/{\rm :}\/

\rritem{(b1)} $b_n(D(Q))=1$ if and only if  $d=n+1$ and up to a linear coordinate change we
have
$\ell_i(x)=x_{i-1}$ for all $i=1,\ldots ,n+1${\rm ;}

\rritem{(b2)} $b_n(D(Q))=2$ if and only if $d=n+2$ and up to a linear coordinate change and
re\/{\rm -}\/ordering of the hyperplanes{\rm ,} $\ell_i(x)=x_{i-1}$ for all
$i=1,\ldots ,n+1$ and
$\ell_{n+2}(x)=x_0+x_1$. 

\endproclaim

Note that the equivalence of (i) and (iii) is a generalization of
Lemma 7 in [Do], and (b1) is a generalization of Theorem 5 in [Do].

To obtain Corollary 4 from Theorem 1 all we need is the following:

\nonumproclaim{Lemma 5}  For any arrangement ${\cal A}$  as above{\rm ,}
$(-1)^n
\chi (D(Q) \setminus H)=b_n(D(Q))$.\endproclaim

Let ${\cal A'}=\{H_i'\}_{i \in I}$ be an affine hyperplane arrangement
in $\C^n$
with complement $M({\cal A'})$    and let $\ell '_i=0$ be an equation
for the hyperplane $H_i'$. Consider the multivalued function
$$ \phi _a:M({\cal A'}) \to \C, ~~~~\phi _a(x)= \prod _i \ell
'_i(x)^{a_i}$$
with $a_i \in \C$. Varchenko conjectured in [V] that for an
essential arrangement ${\cal A'}$ and for generic complex exponents $a_i$ the
function $\phi _a$  has only nondegenerate critical points and their
number is precisely $|\chi(M({\cal A'}))|$. This conjecture was proved
in more general forms by Orlik-Terao [OT3] via algebraic methods and by
Damon [Da2] via topological methods based on [DaM] and [Da1].

In particular Damon shows in Theorem 1 in [Da2] that the function $\phi
_1$ obtained by taking $a_i=1$ for all $i \in I$ has only isolated
singularities and the sum of the corresponding Milnor numbers equals
$|\chi(M({\cal A'}))|$. Consider the morsification
$$\psi (x)= \phi _1(x) - \sum _{j=1}^n b_jx_j$$
where $b_j \in \C$ are generic and small. Then one may think that by
the general property of a morsification, $\psi$ has only nondegenerate
critical points and their number is precisely $|\chi(M({\cal A'}))|$. In
fact, as a look at the simple example $n=3$ and $\phi _1=xyz$ shows,
there are new nondegenerate singularities occurring along the
hyperplanes. This can be restated by saying that
in general one has
$${\rm deg}({\rm grad} \phi _1 )\geq |\chi(M({\cal A'}))|$$
and not an equality similar to our Corollary 4 (1). Note that here
${\rm grad} \phi _1: M({\cal A'}) \to \C ^n$.

The classification of arrangements for which
$|\chi (M({\cal A}'))|=1$ is much more complicated
than the one from Corollary 4(b1) and the interested
reader is referred to [JL].

Theorem 1, in conjunction with Corollary 4, Part (1), has very
interesting consequences.
We say that a topological space $Z$ is {\it minimal} if $Z$ has
the homotopy type of a connected CW-complex $K$ of finite type,
whose number of $k$-cells equals $b_k(K)$ for all $k \in \N$. It is
clear that a minimal space has integral
torsion-free homology. The converse is true for
1-connected spaces; see [PS, Rem. 2.14].

The importance of this notion for the topology of spaces which look
homologically like complements of hyperplane arrangements was
recently noticed in [PS].
Previously, the minimality property was
known only for generic arrangements (Hattori [Hat])
and fiber-type arrangements
(Cohen-Suciu [CS]).
Our next result establishes this property, in full generality. It was
independently  obtained  by Randell [R2,3], using similar techniques.
(See, however,\break Example 13.) The minimality property below should be compared
with the main result from [GM, Part III], where the existence of a
{\it homologically perfect} Morse function is established, for
complements of (arbitrary) arrangements of real affine subspaces;
see [GM, p.\ 236].

\nonumproclaim{{C}orollary 6}  Both complements{\rm ,} $M({\cal A})
\subset \P^n$ and
its cone{\rm ,} $M'({\cal A}) \subset \C^{n+1}${\rm ,} are minimal spaces.\endproclaim

It is easy to see that for $n>1$, the open set $D(f)$ is not minimal
for
$f$ generic of degree $d>1$ (just use
$\pi_1(D(f))=H_1(D(f),\Z)=\Z/d\Z$),
but the Milnor fiber $F$ defined by $f$ is clearly minimal. We do not
know whether the Milnor fiber
$\{ Q=1\}$ associated to an arrangement is   minimal in general.

From Theorem 3 we get a substantial strengthening of some of the main results by Orlik and Terao in [OT2]. Let ${\cal A}'$ be
the central hyperplane arrangement in $\C^{n+1}$ associated to
the projective arrangement ${\cal A}$. Note that $Q(x)=0$ is
a reduced equation for the union $N$ of all the hyperplanes in ${\cal
A}'$.
Let $f \in \C[x_0,\ldots ,x_n]$ be a homogeneous polynomial of degree $e>0$
with
global Milnor fiber $F=\{ x \in \C^{n+1}| f(x)=1\}$ and let
$g: F \setminus N \to \R$
be the function $g(x)=Q(x) {\overline Q}(x)$ associated to
the arrangement. The polynomial $f$ is called
${\cal A}'$-{\it generic} if
\vglue4pt
(GEN1) the restriction of $f$ to any intersection $L$ of
hyperplanes in ${\cal A}'$ is nondegenerate, in the sense that
the associated projective hypersurface in $\P(L)$ is smooth, and
\vglue4pt
(GEN2) the function $g$ is a Morse function.
\vglue4pt
Orlik and Terao have shown in [OT2] that for an {\it essential}
arrangement ${\cal A}'$, the set of ${\cal A}'$-generic functions
$f$ is dense in the set of homogeneous polynomials of degree $e$,
and, as soon as we have an ${\cal A}'$-generic function $f$,
the following basic properties hold for any arrangement.

\ritem{(P1)} $b_q(F,F \cap N)=0$ for $q \not= n$ and

\ritem{(P2)} $b_{n}(F,F \cap N) \leq |C(g)|$, where $C(g)$
is the critical set of the Morse function~$g$.
\vglue4pt 
An explicit formula for the number $|C(g)|$ is given in [OT2] in terms
of the lattice associated to the arrangement ${\cal A}'$.
Moreover, for a special class of arrangements called
{\it pure} arrangements it is shown in [OT2] that (P2) is
actually an equality. In fact, the proof of (P2) in [OT2] uses
Morse theory on noncompact manifolds, but we are unable to see
the details
behind the proof of Corollary (3.5); compare to our discussion in
Example 13.

With this notation  the following is a direct consequence of Theorem 3.

\nonumproclaim{{C}orollary 7} For  any arrangement ${\cal A}$ the
following hold\/{\rm :}\/

\ritem{(i)} the set of ${\cal A}'$\/{\rm -}\/generic functions $f$ is dense
in the set of homogeneous polynomials of degree $e>0${\rm ;}

\ritem{(ii)}  the Milnor fiber $F$ is homotopy equivalent to a {\rm CW}
complex obtained from $F \cap N$ by attaching $|C(g)|$ cells
of dimension $n${\rm ,} where $C(g)$ is the critical set of the Morse
function $g$.
In particular $b_{n}(F,F \cap N)=|C(g)|$ and  the intersection $F \cap
N$ is homotopy equivalent to a bouquet of $|C(g)|-(e-1)^{n+1}$ spheres
$S^{n-1}$.

\endproclaim

Similar results for nonlinear arrangements on complete intersections
have been obtained by Damon in [Da3] where explicit formulas for $|C(g)|$
are given.

The aforementioned results represent a strengthening of those in [D2]
(in which the homological version of Theorems 1 and 3 above was
proved).

The investigation of higher homotopy groups of complements
of complex hypersurfaces (as $\pi_1$-modules) is a very difficult problem.
In the irreducible case, see [Li] for various results on
the first nontrivial higher homotopy group. The arrangements of
hyperplanes provide the simplest nonirreducible situation
(where $\pi_1$ is never trivial, but at the same time rather well
understood). This is the topic of the second part of our paper.

Our results here use the general approach by minimality from [PS],
and significantly extend the homotopy
computations therefrom.
In Section 5, we present a unifying framework for all known explicit
descriptions of nontrivial higher homotopy groups of arrangement
complements,
together with a numerical $K(\pi, 1)$-test. We give specific examples,
in Section 6, with emphasis on combinatorial determination. A general
survey of Sections 5 and 6 follows. (To avoid overloading the
exposition, formulas will be systematically skipped.)

Our first main result in Sections 5 and 6 is Theorem 16. It applies to arrangements
${\cal A}$ which are $k$-{\it generic sections}, $k \geq 2$, of
aspherical arrangements, $\widehat{\cal A}$.
Here '$k$-generic' means, roughly
speaking, that ${\cal A}$ and $\widehat{\cal A}$ have the same
intersection
lattice, up to rank $k+1$; see Section 5(1) for the precise definition.

The {\it general position} arrangements from [Hat] and the {\it
fiber-type}
aspherical ones from [FR] belong to the {\it hypersolvable} class from
[JP1].
Consequently ([JP2]), they all are $2$-generic sections of fiber-type
arrangements. At the same time, the {\it iterated generic hyperplane sections},
${\cal A}$,
of essential aspherical arrangements, $\widehat{\cal A}$, from [R1], are
also
particular cases of $k$-generic sections, with $k={\rm rank}({\cal A})-1$.

For such a $k$-generic section ${\cal A}$, Theorem 16 firstly says
that the
complement $M({\cal A})$ ($M'({\cal A})$) is aspherical if and only if
$p=\infty$,
where $p$ is a topological invariant introduced in [PS]. Secondly (if
$p<\infty$),
one can write down a\break $\Z\pi_1$-module presentation for $\pi_p$, the
first
higher nontrivial homotopy group of the complement (see \S 5(8), (9) for
details). Both results essentially follow from Propositions 14 and 15,
which together imply that $M({\cal A})$ and $M(\widehat{\cal A})$
share the same $p$-skeleton.

In Theorem 18, we substantially
extend and improve
results from [Hat] and [R1] (see also Remark 19). Here we examine
${\cal A}$,
an iterated generic hyperplane section of rank $\geq 3$, of an
essential
aspherical arrangement, $\widehat{\cal A}$. Set $M=M({\cal A})$. In this
case,
$p={\rm rank}({\cal A})-1$ [PS]. We show that the
$\Z\pi_1(M)$-
presentation of $\pi_p(M)$ from Theorem 16 extends to a finite,
minimal,
free $\Z\pi_1(M)$-resolution. We infer that $\pi_p(M)$ cannot be a
projective
$\Z\pi_1(M)$-module, unless ${\rm rank}({\cal A}) = {\rm rank}(\widehat{\cal A})
-1$,
when it is actually $\Z\pi_1(M)$-free.

In Theorem 18 (v), we go beyond the first nontrivial higher homotopy
group.
We obtain a complete description of {\it all} higher rational homotopy
groups,
$L_* := \mathbold{\oplus}_{q \geq 1} ~~\pi_{q+1}(M) \otimes \Q$,
including both the graded Lie algebra structure of $L_*$ induced by the
Whitehead product, and the graded $\Q \pi_1(M)$-module structure.

The computational difficulties related to the twisted homology
of a connected CW-complex (in particular, to the first nonzero higher
homotopy group) stem from the fact that the $\Z \pi_1$-chain complex
of the universal cover is very difficult to describe, in general. As
explained in the introduction, we have two results in this direction,
at the $I$-adic associated graded level: Theorem 20 and Corollary 21.

Corollary 21 belongs to a recurrent theme of our paper:
exploration of  new phenomena of combinatorial determination
in the homotopy theory of arrangements. Our combinatorial determination
property from Corollary 21 should be compared with a fundamental result
of Kohno [K], which says that the rational graded Lie algebra associated
to the lower central series filtration of $\pi_1$ of a
projective hypersurface complement is determined by
the cohomology ring.

In Theorem 23, we examine the hypersolvable arrangements
for which $p = {\rm rank} ({\cal A}) -1$. We establish
the combinatorial determination property of
the $I$-adic associated graded module (over $\Z$) of the first higher
nontrivial homotopy group of the complement, $\pi_p$, in Theorem 23 (i).
The proof uses in an essential way the ubiquitous Koszul property
from homological algebra.

We also infer from Koszulness, in Theorem 23 (ii), that the successive
quotients of the $I$-adic filtration on $\pi_p$ are finitely generated free
abelian groups, with ranks given by the combinatorial
$I$-{\it adic filtration formula} $(22)$. This resembles the
{\it lower central series} ({\rm LCS}) {\it formula}, which expresses
the ranks of the quotients of the lower central series of $\pi_1$
of certain arrangements, in combinatorial terms. The {\rm LCS} formula
for pure braid groups was discovered by Kohno,
starting from his pioneering work in [K]. It was established for all
fiber-type arrangements in [FR],
and then extended to the hypersolvable class in [JP1].

Another new example of combinatorial determination is the fact that
the generic affine part of a union of hyperplanes has the homotopy type of
the Folkman complex, associated to the intersection lattice.
This follows from Theorem 1 and Corollary 4; see the discussion after the
proof of Theorem 3.

\vglue-6pt
\section{ Polar curves, affine Lefschetz theory\\ and degree of gradient
maps}

 \vglue-6pt
 
The use of the local polar varieties in the study of singular spaces
is already a classical subject;  see L\^e [Le1],  L\^e-Teissier [LT]
and the references therein.
Global polar curves in the study of the topology of polynomials
is a topic
under intense investigations; see for instance Cassou-Nogu\`es and
Dimca [CD], Hamm [H], N\'emethi [N1,2], Siersma and Tib\u ar [ST].
For all the proofs in this paper, the classical theory is sufficient:
indeed, all the objects being homogeneous, one can localize at the
origin of
$\C^{n+1}$ in the standard way, see [D1]. However, using
geometric
intuition, we find it easier to work with global objects,
and hence we adopt this viewpoint in the sequel.

We recall briefly the notation and the results from [CD], [N1,2].
Let $h \in \C[x_0,\ldots , x_n]$ be a polynomial (even nonhomogeneous
to start with) and assume that the
fiber $F_t=h^{-1}(t)$ is smooth, for some fixed $t \in \C$.

 For any hyperplane in $\P^n$, $ H:\ell =0$ where
 $  \ell(x)=c_0x_0+c_1x_1+\cdots +c_nx_n $,
we define the corresponding polar variety $\Gamma _H $ to be the union
of the irreducible components of the variety
$$ \{x \in \C^{n+1} \ |   \  {\rm rank}(dh(x),d\ell (x))=1 \} $$
which are not contained in the critical set $S(h) = \{x \in \C^{n+1} \
| \
dh(x)=0 \} $ of~$h$.

\nonumproclaim{Lemma 8 {\rm (see [CD], [ST])}}  For a generic hyperplane
$H${\rm ,}

\ritem{(i)} The polar variety $\Gamma _H$ is either empty or a curve\/{\rm ;}
i.e.{\rm ,} each
irreducible component of $\Gamma _H$ has dimension $1$.

\ritem{(ii)} {\rm dim}$(F_t \cap \Gamma _H) \leq 0$ and the intersection
multiplicity
$(F_t,\Gamma_H)$ is independent of $H$.

\ritem{(iii)} The multiplicity $(F_t, \Gamma _H)$ is equal to the number of
tangent hyperplanes to $F_t$ parallel to the hyperplane $H$. For each
such tangent hyperplane $H_a${\rm ,} the intersection $ F_t \cap H_a$ has
precisely one singularity{\rm ,} which is an ordinary double point.

\endproclaim

 The nonnegative integer $(F_t, \Gamma _H)$ is called the polar
invariant
of the hypersurface $F_t$
and is denoted by $P(F_t)$. Note that $P(F_t)$ corresponds exactly
to the classical notion of
class of a projective hypersurface; see [L].

We think of a projective hyperplane $H$  as the direction
of an affine hyperplane $H' =\{ x \in \C^{n+1}|  \ell (x)=s \}$ for
$s \in \C$. All the affine hyperplanes with the same direction form a
pencil,
and it is precisely this type of pencil  that is used in the
affine
Lefschetz theory; see [N1,2]. N\'emethi considers only connected affine
varieties, but his results clearly extend to the case of any pure
dimensional smooth variety.

\nonumproclaim{Proposition 9 {\rm(see [CD], [ST])}} For a generic
hyperplane $H'$ in the pencil of all hyperplanes in
$\C^{n+1}$ with a fixed generic direction $H${\rm ,} the fiber $F_t$ is
homotopy
equivalent to a {\rm CW-}\/complex
obtained from  the section $F_t \cap H'$ by
attaching $P(F_t)$  cells of dimension $n$.
In particular 
$$P(F_t)=(-1)^{n}(\chi(F_t)- \chi(F_t \cap H'))=(-1)^n \chi(F_t
\setminus H').$$
\endproclaim

Moreover, in this statement `generic' means that the affine hyperplane
$H'$ has to
verify the following two conditions:

\ritem{(g1)} its direction in $\P^n$ has to be generic, and

\ritem{(g2)} the intersection $ F_t \cap H'$ has to be smooth.
\vglue4pt
These two conditions are not stated in [CD], but the reader should have
no
problem in checking them by using Theorem $3'$ in [CD] and the fact
proved
by N\'emethi in [N1,2] that the only bad sections in a good
pencil are the singular sections. Completely similar results hold for
generic
pencils with respect to a closed smooth subvariety $Y$ in some affine
space
$\C^N$; see [N1,2], but note that the polar curves are not mentioned
there.

\vglue4pt {\it {P}roof of Theorem {\rm 1}}.
In view of Hamm's affine  Lefschetz theory, see\break [H, Th.~5],
the only thing to prove is the equality between the number $k_n$ of
$n$-cells attached and the degree of the gradient.

Assume from now on that the polynomial $h$ is homogeneous of degree
$d$ and that $t=1$. It follows from (g1) and (g2) above that we may
choose the generic hyperplane $H'$ passing through the origin.

Moreover, in this case, the polar curve $\Gamma _H$, being defined by
homogeneous equations, is a union of lines $L_j$ passing through the
origin.
For each such line we choose a parametrization $t \mapsto a_j t$ for
some $a_j \in \C^{n+1}, a_j \not=0$. It is easy to see that the
intersection $F_1 \cap L_j$ is either empty (if $h(a_j)=0$) or consists
of exactly $d$ distinct points with multiplicity one (if
$h(a_j)\not=0$).
The lines of the second type are in bijection with the points in
${\rm grad}(h)^{-1}(D_{H'})$,
where $D_{H'}\in \P^n$ is the point corresponding to the direction of
the hyperplane $H'$. It follows that
$$ d \cdot {\rm deg}({\rm grad}(h))= P(F_1).$$

The $d$-sheeted unramified coverings $ F_1 \to D(h)$ and
$F_1 \cap H' \to D(h) \cap H$ give the result, where $H$ is the
projective hyperplane corresponding to the affine hyperplane (passing
through the origin) $H'$. Indeed, they imply the equalities:
$\chi(F_1)=d \cdot \chi (D(h))$ and
$\chi(F_1 \cap H')=d \cdot \chi (D(h) \cap H)$. Hence we have
${\rm deg}({\rm grad}(h))= (-1)^n \chi (F_1,F_1 \cap H')/d=(-1)^n \chi(D(h),D(h)
\cap H)=k_n.$
\vglue4pt
\pagegoal=50pc

{\it Remark} 10.  The gradient map ${\rm grad}(h)$ has a natural
extension to the
larger open set $D'(h)$ where at least one of the partial derivatives
of
$h$ does not vanish. It is obvious (by a dimension argument) that this
extension has the same degree as the map ${\rm grad}(h)$.

\section{Nonproper Morse theory}
\vglue-6pt

For the convenience of the reader we recall, in the special case  
needed,
a basic result of Hamm, see [H, Prop.~3], with our addition
concerning
the condition (c0) in [DP, Lemma 3 and Ex. 2]. The final
claim on
the number of cells to be attached is also standard, see for instance
[Le1].

\nonumproclaim{Proposition 11}  Let $A$ be a smooth
algebraic subvariety in $\C^p$ with ${\rm dim}A=m$.
Let $f_1,\ldots ,f_p $ be polynomials in $\C[x_1,\ldots ,x_p]$.
For $ 1 \leq j \leq p${\rm ,} denote by  $\Sigma _j$ the set of critical
points
of the mapping $(f_1,\ldots ,f_j): A \setminus \{z \in A \ | \ f_1(z)=0 \} \to
\C^j$
and let $\Sigma _j'$ denote the closure of $\Sigma _j$ in $A$.
Assume  that the following conditions hold.
\vglue2pt
\item{\rm (c0)} The set $\{z \in A \ | \ |f_1(z)| \leq a_1,\ldots ,  |f_p(z)| \leq a_p \}$
is compact for any positive numbers $a_j$, $j=1,\ldots ,p$.
\vglue2pt

\item{\rm (c1)} The  critical set $\Sigma _1$ is finite.

\vglue2pt
\item{\rm (cj)} {\rm (for $j=2,\ldots ,p$)} The map  $(f_1,\ldots ,f_{j-1}): \Sigma _j' \to \C
^{j-1}$
is  proper.
\vglue2pt

Then $A$ has the homotopy type of a space obtained from
$A_1=\{z \in A \ | \ f_1(z)=0\}$ by attaching $m$-cells and the number
of these cells is the sum of the Milnor numbers $\mu(f_1,z)$ for
$z \in \Sigma _1$. 
\endproclaim

\pagegoal=50pc
 
{\it Proof of Theorem {\rm 3}}.
We set $X=h^{-1}(1)$.
Let  $v: \C^{n+1} \to \C^N$ be the Veronese mapping of degree $e$
sending $x$ to all the monomials of degree $e$ in $x$ and set $Y=v(X)$.
Then $Y$ is a smooth closed subvariety in $\C^N$ and $v:X \to Y$ is an
unramified  covering of degree $c$, where $c={\rm g.c.d.}(d,e)$.
To see this, use the fact that $v$ is a closed immersion on
$\C^N \setminus \{0\}$ and $v(x)=v(x')$ if and only if
$x'=u \cdot x$ with
$u^c=1$.

Let $H$ be a generic hyperplane direction in $\C^N$ with respect to the
subvariety $Y$ and let $C(H)$ be  the finite set of all the points
$p \in Y$ such that there is an affine hyperplane $H'_p$ in the
pencil determined by $H$ that is tangent to $Y$ at the point $p$ and
the intersection $Y \cap H'_p$ has a complex Morse
(alias nondegenerate, alias $A_1$) singularity.
Under the Veronese mapping $v$, the generic hyperplane direction
$H$ corresponds to a homogeneous polynomial of degree $e$ which
we call from now on $f$.

To prove the first claim (i) we proceed as follows. It is known that
using affine Lefschetz theory for a pencil of hypersurfaces $\{h=t\}$
is equivalent to using (nonproper) Morse theory for the function $|h|$
or,
what amounts to the same, for the function $|h|^2$. More explicitly,
in view of the last statement at the end of the proof of Lemma (2.5) in
[OT2] (which clearly applies to our more general setting since all
the computations there are local), $g$ is a Morse function
if and only if each
critical point of $h:F \setminus N \to \C$ is an $A_1$-singularity.
By the homogeneity of both $f$ and $h$, this last condition on $h$
is equivalent to the fact that each critical point of the function
$f:X \to \C$ is an $A_1$-singularity, condition fulfilled in view of
the choice of $H$ and since $v:X \to Y$ is a local isomorphism.

Now we pass on to the proof of the claim (ii) in Theorem 3. One can
derive this claim easily from Proposition 9.14 in Damon [Da1].
However since his proof is using
previous results by Siersma [Si] and Looijenga [Lo], we think that our
original proof given below and based on Proposition 11, longer but more
self-contained, retains its interest.

Any polynomial function $h:\C^{n+1} \to \C$ admits
a Whitney stratification satisfying the Thom $a_h$-condition. This is
a constructible stratification ${\cal S}$ such that the open stratum,
say $S_0$, coincides with the set of regular points for $h$ and for
any other stratum, say $S_1 \subset h^{-1}(0)$, and any
sequence of points $q_m \in S_0$ converging to $q \in S_1$
such that the sequence of tangent spaces $T_{q_m}(h)$ has a limit
$T$; one has $T_qS_1 \subset T$. See Hironaka [Hi, Cor.~1, p.~248].
The requirement of $f$ proper in that corollary is not
necessary in our case, as any algebraic map can be compactified. Here
and in the sequel, for a map $\phi : {\cal X} \to {\cal Y}$ and a point
$q \in {\cal X}$ we denote by $T_q(\phi)$ the tangent space to the
fiber
$\phi ^{-1}(\phi (q))$ at the point~$q$, assumed to be a smooth point
on
this fiber.

Since in our case $h$ is a homogeneous polynomial, we can find a
stratification ${\cal S}$ as above such that all of its strata are
$\C^*$-invariant, with respect to the natural $\C^*$-action on
$\C^{n+1}$.
In this way we obtain an induced Whitney stratification  ${\cal S'}$ on
the
projective hypersurface $V(h)$.
We select our polynomial $f$ such that the corresponding projective
hypersurface $V(f)$ is smooth and transversal to the stratification
${\cal S'}$. In this way we get  an induced Whitney stratification
${\cal S'}_1$ on the projective complete intersection $V_1= V(h) \cap
V(f)$.

We use Proposition 11 above with $A=F$ and $f_1=h$. All we have to
show
is the existence of polynomials $f_2,\ldots ,f_{n+1}$ satisfying the
conditions
listed in Proposition 11.

We will select these polynomials inductively to be generic linear forms
as follows. We choose $f_2$ such that the corresponding hyperplane
$H_2$
is transversal to the stratification  ${\cal S'}_1$. Let ${\cal S'}_2$
denote the induced stratification on $V_2 =V_1 \cap H_2$. Assume that
we have constructed $f_2,\ldots ,f_{j-1}$, ${\cal S'}_1,\ldots ,{\cal
S'}_{j-1}$
and $V_1,\ldots ,V_{j-1}$. We choose $f_j$ such that the corresponding
hyperplane
$H_j$ is transversal to the stratification  ${\cal S'}_{j-1}$. Let
${\cal S'}_j$ denote the induced stratification on $V_j =V_{j-1} \cap
H_j$.
Do this for $j=3,\ldots ,n$ and choose for $f_{n+1}$ any linear form.

With this choice it is clear that for $1 \leq j \leq n$, $V_j$ is a
complete intersection of dimension $n-1-j$. In particular,
$V_n=\emptyset$;
i.e.
$$  \{x \in \C^{n+1} \ | \
f(x)=h(x)=f_2(x)=\cdots = f_{n}(x)=0\}=\{0\}~.\speqnu{\rm  c0'}$$

Then the map $(f,h,f_2,\ldots ,f_n):\C^{n+1} \to \C^{n+1}$ is proper,
which clearly implies the condition (c0).
 
The condition (c1) is fulfilled by our construction of $f$.
Assume that we have already checked that the conditions $(ck)$
are fulfilled for $k=1,\ldots ,j-1$. We explain now why  the next
condition (cj) is fulfilled.

Assume that the condition (cj) fails. This is equivalent to
the existence of a sequence $p_m$ of points in $\Sigma _j'$ such that

\ritem{$(*)$} $|p_m| \to \infty$ and $f_k(p_m) \to b_k$ (finite limits) for
$1 \leq k \leq j-1$.
\vglue4pt 
Since $\Sigma _j$ is dense in
$\Sigma _j'$, we can even assume that $p_m \in  \Sigma _j$.

Note that $\Sigma _{j-1} \subset \Sigma _j$ and the condition
$c(j-1)$ is fulfilled. This implies that we may choose our sequence
$p_m$ in the difference
$\Sigma _j \setminus  \Sigma _{j-1}$. In this case we get
\ritem{$(**)$} $f_j \in {\rm Span}(df(p_m),dh(p_m),f_2,\ldots ,f_{j-1})$
\vglue4pt\noindent the latter being a $j$-dimensional vector space.

Let $q_m ={p_m \over |p_m|} \in S^{2n+1}$. Since the sphere
$S^{2n+1}$ is compact we can assume that the sequence $q_m$
converges to a limit point $q$. By passing to the limit in $(*)$
we get $q \in V_{j-1}$. Moreover, we can assume (by passing to a
subsequence) that the sequence of $(n-j+1)$-planes
$T_{q_m}(h,f,f_2,\ldots ,f_{j-1})$ has a limit $T$. Since
$p_m \notin \Sigma_{j-1}$, we have
$$ T_{q_m}(h,f,f_2,\ldots ,f_{j-1})=T_{q_m}(h) \cap T_{q_m}(f) \cap H_2
\cap\cdots \cap H_{j-1} ~.$$
As above, we can assume that the sequence $T_{q_m}(h)$ has a limit
$T_1$ and, using the $a_h$-condition for the stratification ${\cal S}$
we get $T_qS_i \subset T_1$ if $q \in S_i$. Note that
$T_{q_m}(f) \to T_{q}(f)$ and hence
$T=T_1 \cap T_{q}(f) \cap H_2 \cap\cdots \cap H_{j-1}$. It follows that
$$T_qS_{i,j-1}= T_qS_i \cap T_{q}(f) \cap H_2 \cap\cdots \cap H_{j-1}
\subset T ~,$$
where $S_{i,j-1}= S_i \cap V(f) \cap H_2 \cap\cdots \cap H_{j-1}$ is
the stratum corresponding to the stratum $S_i$ in the stratification
${\cal S'}_{j-1}$. On the other hand, the condition $(**)$ implies that
$T_qS_{i,j-1} \subset T \subset   H_j $, a contradiction to the fact
that
$H_j$ is transversal to ${\cal S'}_{j-1}$.

To prove (iii) just note that the intersection $F \cap N$ is
$(n-2)$-connected (use the exact homotopy sequence of the pair $(F,F \cap N)$
and the fact that $F$ is $(n-1)$-connected) and has the homotopy type of
a CW-complex of dimension $\leq (n-1)$.
\enddemo

Let us now reformulate slightly  Theorem 1 as explained already in
Section~1. Note that $\chi (D(h) \setminus H)=\chi (\P^n \setminus (V(h)
\cup H))= \chi (\C ^n \setminus V(h)_a)=
1-(1+(-1)^{n-1}b_{n-1}(V(h)_a))=(-1)^n b_{n-1}(V(h)_a)$.
Here $V(h)_a=V(h) \setminus H$ is the affine part of the projective
hypersurface $V(h)$ with respect to the hyperplane at infinity $H$ and,
according to L\^e's affine Lefschetz Theorem, is homotopy equivalent to a
bouquet of $(n-1)$-spheres. It follows that
$$ {\rm deg}({\rm grad}(h))=b_{n-1}(V(h)_a) ~.$$
Even in the case of an arrangement ${\cal A }$,
the corresponding equality, $b_n(M({\cal A}))=b_{n-1}({\cal A}_a)$,
derived by using Corollary 4, seems to be new. Here we denote by ${\cal A }$
not only the projective arrangement but also the union of all the
hyperplanes in ${\cal A }$. Theorem 4.109 in [OT1] implies that, in the case
of an essential arrangement, the bouquet of spheres ${\cal A }_a$
considered above is homotopy equivalent to the Folkman complex $F({\cal A
}')$ associated to the corresponding central arrangement
 ${\cal A }'$  in $\C^{n+1}$; see [OT1, pp.~137--142] and [Da1, p.~40].

The main interest in Dolgachev's paper [Do] is focused on homaloidal
polynomials, i.e., homogeneous polynomials $h$ such that $
{\rm deg}({\rm grad}(h))=1$. In view of the above reformulation of Theorem 1 it follows that a
polynomial $h$ is homaloidal
if and only if the affine hypersurface $V(h)_a$ is homotopy equivalent
to an $(n-1)$-sphere. There are several direct consequences of this
fact.
\vglue4pt
(i) If $V(h)$ is either a smooth quadric (i.e. $d=2$) or the union of a
smooth
quadric $Q$ and a tangent hyperplane $H_0$ to $Q$ ( in this case
$d=3$), then
$ {\rm deg}({\rm grad}(h))=1.$ Indeed, the first case is obvious (either from the
topology or the algebra), and the second case follows from the fact that
both $H_0 \setminus H$ and $(Q \cap H_0) \setminus H$ are contractible.
\vglue4pt
(ii) Using the topological description of an irreducible projective
curve as the wedge of a smooth curve of genus $g$ and $m$ circles, we see
that for an irreducible plane curve $C$ of degree $d$,
$$b_1(C_a)=2g+m+d-1.$$
Using this and the Mayer-Vietoris exact sequence for homology one can
derive   Theorem~4 in [Do], which gives the list of all reduced
homaloidal polynomials $h$  in the case $n=2$. This list is reduced to the
two examples in (i) above plus the union of three nonconcurrent lines.
\vglue4pt
(iii) If the hypersurface $V(h)$ has only isolated singular points, say
at  $a_1,\ldots ,a_m$, then our formula above gives
$${\rm deg}({\rm grad}(h))=(d-1)^n- \sum_{i=1}^m \mu (V(h),a_i) ~,$$
where $\mu (V,a)$ is the Milnor number of the isolated hypersurface
singularity $(V,a)$; see [D1, p.~161]. We conjecture that when $n>2$ and
$d>2$ one has
${\rm deg}({\rm grad}(h))>1$ in this situation, unless $V(h)$ is a cone and then of
course ${\rm deg}({\rm grad}(h))=0$
.
For more details on this conjecture and its relation to the work by A. duPlessis and C. T. C. Wall in  [dPW], see  [D3].

\section{Complements of hyperplane arrangements}

{\it Proof of Lemma} 5.
We are going to derive this easy result from\break [H, Th.~5], by using the
key homological features of arrangement complements. By Hamm's theorem,
the equality between $(-1)^n \chi(D(Q)\setminus H)=(-1)^n \chi(D(Q)$,\break
$D(Q)\cap H)$
and $b_n(D(Q))$ is equivalent to
$$ b_{n-1}(D(Q)\cap H)=b_{n-1}(D(Q)).\eqno(*)$$
All we can say in general is that $b_{n-1}(D(Q)\cap H) \geq
b_{n-1}(D(Q))$.
In the arrangement case, the other inequality follows from two standard
facts
(see [OT1, Cor.~5.88 and Th.~5.89]): $H^{n-1}(D(Q)\cap H)$ is
generated by products of cohomology classes of degrees $< n-1$ (if
$n>2$);
$H^1(D(Q)) \to H^1(D(Q)\cap H)$ is surjective (which in particular
settles
the case $n=2$). See also Proposition~2.1 in [Da3].

\vglue8pt {\it Proof of Corollary {\rm 4}}.
 (1) This claim follows directly from Theorem 1 and Lemma 5.
\vglue4pt 
(2)(a) To complete this proof  we only have to explain why the claims
(ii) and (iii) are equivalent. This in turn is an immediate
consequence
of the well-known equality: ${\rm deg} P_{\cal A}(t)={\rm codim} (\mathbold{\cap}_{i=1}^d
H_i)-1$,
where $P_{\cal A}(t)$ is the Poincar\' e polynomial of $D(Q)$; see
[OT1,
Cor.~3.58, Ths.~3.68 and 2.47].
\vglue4pt
(2)(b) The inequality can be proved by induction on $d$ by the
method of deletion and restriction; see [OT1, p.~17].
\vglue6pt

{\it  Proof of Corollary} 6.
Using the affine Lefschetz theorem of Hamm  (see Theorem 5 in [H]),
we know that for a generic projective hyperplane $H$, the space $M$ has
the homotopy type of a space obtained from $M \cap H$ by attaching
$n$-cells.
The number of these cells is given by
\vglue3pt
\centerline{$(-1)^n \chi (M, M \cap H)=(-1)^n \chi (M \setminus H)=b_n(M);$}
\vglue3pt\noindent 
see Corollary 4 above.

To finish the proof of the minimality of $M$ we proceed by induction.
Start with a minimal cell structure, $K$, for $M\cap H$, to get a cell
structure, $L$, for $M$, by attaching $b_n(M)$ top cells to $K$. By
minimality, we know that $K$ has trivial cellular incidences. The fact
that
the number of top cells of $L$ equals $b_n(L)$ means that these cells
are
attached with trivial incidences, too, whence the minimality of $M$.

Finally, $M'$ has the homotopy type of $M\times S^1$, being therefore
minimal, too.
 
\vfil {\it {R}emark {\rm 12}}. (i) Let $\mu_e$ be the cyclic group of the
$e$-roots of unity.
Then there are natural algebraic actions of $\mu_e$ on the spaces
$F \setminus N$  and $F \cap N$ occurring in Theorem 3. The
corresponding weight
equivariant Euler polynomials (see [DL] for a definition) give
information on the relation between the induced\break $\mu_e$-actions
on the cohomology $H^*(F \setminus N,\Q)$ and $H^{n-1}(F \cap N,\Q)$
and the functorial Deligne
mixed Hodge structure present on cohomology.

When $N$ is a hyperplane arrangement ${\cal A}'$  and $f$ is an
${\cal A}'$-generic function, these weight equivariant Euler
polynomials
can be combinatorially computed from the lattice associated to the
arrangement
(see Corollary (2.3) and\break Remark (2.7) in [DL]) by the fact that the
 weight equivariant Euler polynomial of the $\mu_e$-variety $F$ is
known;
 see  [MO] and [St]. This gives in particular the characteristic
polynomial
of the monodromy associated to the function\break $f: N \to \C$.
\eject

(ii) The  minimality property turns out to be useful
in the context of {\it homology with twisted coefficients}. Here is
a simple example. (More results along this line will be published
elsewhere.)
Let $X$ be a connected CW-complex of finite type. Set $\pi :=
\pi_1(X)$.
For a left $\Z \pi$-module $N$, denote by $H_*(X, N)$ the homology of
$X$ with
local coefficients corresponding to $N$. One knows that  $H_*(X, N)$
may be computed as
the homology of the chain complex $C_*(\widetilde X)\otimes_{\Z \pi}
N$, where
$C_*(\widetilde X)$ denotes the $\pi$-equivariant chain complex of the
universal cover of $X$; see [W, Ch. VI].

Assume now that $X$ is minimal. If $N$ is a finite-dimensional
$\K$-represen\-tation of $\pi$
over a field $\K$, we obtain, from the above description of twisted
homology, that
$${\rm dim}_{\K}~H_q(X, N)~\leq~({\rm dim}_{\K}~N)\cdot b_q(X)~,~\hbox{ for all } q~.$$
\noindent When $X$ is, up to homotopy, an arrangement complement, and
$\K = \C$, we thus recover the main result of [C]. (Twisted cohomology
may be treated similarly.)
\vglue8pt
 
 {\it Example} 13.  In this example we explain why special care
is needed when doing Morse theory on noncompact manifolds as in [OT2]
and [R2].
Let's start with a very simple case, where computations are easy.
Consider the Milnor fiber,
$X \subset \C^2$, given by $\{ xy=1 \}$. The
hyperplane $\{ x+y=0 \}$ is generic with respect to the arrangement
$\{ xy=0 \}$, in the sense of [R2, Prop.~2]. Set $\sigma := |x+y|^2$.
When trying to do proper Morse theory with boundary, as in [R2, Th.~3],
one faces a delicate problem on the boundary. Denoting by $B_R$ ($S_R$)
the closed ball (sphere) of radius $R$ in $\C^2$, we see easily
that
$X \cap B_R = \emptyset$, if $R^2<2$, and that the intersection $X \cap
S_R$
is not transverse, if $R^2=2$. In the remaining case ($R^2>2$), it is
equally
easy to check that the restriction of $\sigma$ to the boundary $X\cap
S_R$ always
has  eight critical points (with $\sigma \not= 0$). In our very
simple
example, all these critical points are `\` a gradient sortant' (in the
terminology of [HL, Def.~3.1.2]). The proof of this fact does not seem
obvious, in general. At the same time, this property seems to be
needed, in
order to get the conclusion of [R2, Th.~3] (see [HL, Th.~3.1.7]).

One can avoid this problem as follows.
(See also Theorem 3 in [R3].)
Let $X$ denote the affine Milnor fiber and
$\ell :X  \to \C$ the linear function induced by the equation of any
hyperplane. Then for a real number $r>0$, let $D_r$
be the open disc $|z|<r$ in $\C$.
The function $\ell$ having a finite number
of critical points, it follows that $X$ has
the same homotopy type (even diffeo type)
as the cylinder $X_r=\ell^{-1}(D_r)$ for
$r>>0$. Fix such an $r$. For $R>>r$, we have that
$Y=X_r \cap B_R$ has the same homotopy
type as $ X_r$.

Moreover, if the hyperplane $\ell =0 $ is generic
in the sense of L\^e/N\'emethi, it follows
that the real function $ |\ell |$ has no critical
points on the boundary of $Y$ except those
corresponding to the minimal value 0 which do not matter.\break  
In  the example treated before, one can check
that the critical values\break corresponding to the eight
critical points tend to infinity when $R \to \infty$;
i.e., the eight singularities are no longer in $Y$ for a good choice of $r$
and $R$ !

Note that $Y$ is
 a noncompact manifold with a noncompact boundary, but the sets
$Y^{r_0}=\{ y \in Y| |\ell(y)| \leq r_0 \}$ are
compact for all $0 \leq r_0 <r$.
If $R$ is chosen large enough, then $Y^0$ has the same homotopy type as
$X \cap \{\ell =0\}$, and hence we can use proper
Morse theory on manifolds with boundary to get the result.

The idea behind Proposition 11 is the same: one can do
Morse theory on a noncompact manifold with corners
if there are no critical points on the boundary,
and this is achieved by an argument similar
to the above and based on the conditions (cj); see [H] for more
details.

\section{$k$-generic sections of aspherical arrangements}

The preceding minimality result (Corollary 6) enables us to use the
general
method of [PS] to get explicit information on higher homotopy groups of
arrangement complements, in certain situations.
We begin by describing a framework that encompasses
all such known computations. (For the basic facts in arrangement
theory, we
use reference [OT1].)

Let ${\cal A}=\{ H_1,\dots ,H_n \}$ be a projective hyperplane
arrangement in
$\P (V)$, with associated central arrangement, ${\cal A}'=\{ H_1',\dots
,H_n'\}$,
in $V$. Let $M({\cal A}) \subset \P (V)$ and $M'({\cal A}) \subset V$
be the
corresponding arrangement complements.
Denote by ${\cal L}({\cal A})$ the intersection lattice, that is the
set of intersections of hyperplanes from ${\cal A}'$, $X$ (called {\it
flats}),
ordered by reverse inclusion. We will assume that ${\cal A}$ is
essential.
(We may do this, without changing the homotopy types of the complements
and the
intersection lattice, in a standard way; see [OT1, p.~197].) Then there
is a canonical Whitney stratification of $\P (V)$, ${\cal S}_{\cal A}$,
whose
strata are indexed by the nonzero flats, having $M({\cal A})$ as top
stratum;
see [GM, III 3.1 and III 4.5]. Thus, the genericity condition from
Morse theory takes a particularly simple form, in the arrangement case.

To be more precise, we will need the following definition.
Let $U\subset V$ be a complex vector subspace. We say that $U$ is
${\cal L}_k
({\cal A})$-{\it generic} ($0\leq k < {\rm rank}({\cal A})$) if
$$ {\rm codim}_V(X)={\rm codim}_U(X\cap U),~ \hbox{ for all } X \in {\cal L}({\cal
A}),
~{\rm codim}_V(X)\leq k+1.\eqno(1)$$
 It is not difficult to see that $(1)$ forces
$k\leq {\rm dim}\ U-1$ and that
$\P (U)$ is transverse to ${\cal S}_{\cal A}$ if and only if $U$ is
${\cal L}_k({\cal A})$-generic, with $k={\rm dim}\ U-1$. Consider also the
restriction,
${\cal A}^U := \{ \P (U)\cap H_1,\dots , \P (U)\cap H_n \}$, with
complement $M({\cal A})\cap \P (U)$. A direct application of [GM, Th.~II 5.2]
gives then the following:

\ritem{(2)} If $k={\rm dim}\ U-1$ and $U$ is ${\cal L}_k({\cal A})$-generic, then
the pair
$(M({\cal A}),\break M({\cal A})\cap \P (U))$ is $k$-connected.
\eject

The next proposition, which will provide our
framework for homotopy computations by minimality, upgrades the above
implication $(2)$
to the level of cells, for the case of an arbitrary ${\cal L}_k({\cal
A})$-generic
section, $U$. 

\nonumproclaim{Proposition 14} Let ${\cal A}$ be an arrangement
in $\P (V)${\rm ,}
and let $U\subset V$ be a subspace. Assume that both ${\cal A}$ and the
restriction
${\cal A}^U$ are essential.
If $U$ is ${\cal L}_k({\cal A})$\/{\rm -}\/generic
$(0\leq k < {\rm rank}({\cal A})),$ then the inclusion,
$ M({\cal A})\cap \P (U)\subset M({\cal A})${\rm ,}
has the homotopy type of a cellular map, $j:X \to Y${\rm ,} where\/{\rm :}\/
\ritem{(i)} Both $X$ and $Y$ are minimal {\rm CW-}\/complexes.

\ritem{(ii)} At the level of $k$\/{\rm -}\/skeletons{\rm ,} $X^{(k)}=Y^{(k)}${\rm ,} and the
restriction of
$j$ to $X^{(k)}$ is the identity.

\endproclaim

{\it Proof.} If $k< {\rm dim}\ U-1$, we may find a hyperplane $H\subset V$
which
is  ${\cal L}_m({\cal A})$-generic, $m={\rm dim}\ H-1$, and whose trace on $U$,
$H\cap U$, is ${\cal L}_{m'}({\cal A}^U)$-generic, $m'={\rm dim}~(H\cap
U)-1$. It follows that
${\cal A}^H$ and $H\cap U$ satisfy the assumptions of the proposition,
with the same
$k$. The minimal CW-structures for\break $M({\cal A})\cap \P (H\cap U)$ and
$M({\cal A})\cap \P (H)$ may be extended to minimal structures for
$M({\cal A})\cap \P (U)$
and $M({\cal A})$, exactly as in the proof of Corollary 6. The second
claim of the
proposition follows again by induction, together with routine
homotopy-theoretic arguments,
involving cofibrations and cellular approximations.

If  $k = {\rm dim}\ U-1$, it is not difficult to see that ${\cal A}^U$ may be
obtained from
${\cal A}$ by taking ${\rm dim}\ V - {\rm dim}\ U$ successive generic hyperplane
sections.
Therefore, the method of proof of Corollary 6 shows that, up to
homotopy,
$M({\cal A})\cap \P (U)$ is the $k$-skeleton of a minimal
CW-structure for $M({\cal A})$.
\vglue12pt

Let ${\cal A}$ be an arrangement in $\P (V)$, and let $U\subset V$ be a
proper subspace
which is ${\cal L}_0({\cal A})$-generic. Assume that both ${\cal A}$
and ${\cal A}^U$
are essential. The preceding proposition leads to the following
combinatorial definition:

$$k({\cal A}, U):= {\rm sup}~ \{ 0\leq \ell < {\rm rank} ({\cal A})~\mid~
U~~{\rm is}~~
{\cal L}_{\ell}({\cal A})\hbox{-{\rm generic}} \}~,\eqno(3) $$
and to the next topological counterpart:
$$ p({\cal A}, U):= {\rm sup}  \{q \geq 0~\mid~ b_r(M({\cal A})) =
b_r(M({\cal A}^U))~,~ \hbox{ for all } r\leq q \}~.\eqno(4)$$
 
\nonumproclaim{Proposition 15} $k({\cal A}, U) = p({\cal A}, U)$.
\endproclaim

 {\it Proof.} The inequality  $k({\cal A}, U)\leq p({\cal A}, U)$
follows from
Proposition 14. Denoting by $r$ and $r'$ the ranks of ${\cal A}$ and
${\cal A}^U$
respectively, we know that  $b_s(M({\cal A}^U))=0$, for $s\geq r'$, and
$b_s(M({\cal A}))\neq 0$, for $s<r$. It follows that  $p({\cal A},
U)<r'$,
since $r'<r$. To show that  $p({\cal A}, U)\leq k({\cal A}, U)$, we
have to
pick an arbitrary independent subarrangement, ${\cal B}\subset {\cal
A}$,
of $q+1$ hyperplanes, $1\leq q\leq  p({\cal A}, U)$, and verify that
the restriction
${\cal B}^U$ is independent. To this end, we will use three well-known
facts
(see [OT1]). Firstly, the natural map, $H^*M({\cal A})\to H^*M({\cal
A}^U)$,
is onto. Secondly, the natural map, $H^*M({\cal B})\to H^*M({\cal A})$,
is monic. Together with definition $(4)$ above, these two facts imply
that the natural map,
$H^*M({\cal B})\to H^*M({\cal B}^U)$, is an isomorphism, up to degree
$q$.
The combinatorial description of the cohomology algebras of arrangement
complements by Orlik-Solomon algebras may now be used to deduce the
independence of
${\cal B}^U$ from   ${\cal B}$.

\vglue12pt

\centerline{\bf Description of the $k$-generic framework}

\vglue8pt

We want to apply Theorem 2.10 from [PS] to an (essential) arrangement
${\cal A}$
in $\P (U)$, with complement $M:=M({\cal A})$.
Set $\pi=\pi_1 (M)$. The method of [PS] requires both $M$
and $K(\pi, 1)$ to be minimal spaces, with cohomology algebras
generated
in degree $1$. By Corollary 6 (and standard facts in arrangement
cohomology) all these assumptions will be satisfied, as soon as $K(\pi,
1)$ is
(up to homotopy) an arrangement complement, too.

This in turn happens whenever ${\cal A}$ is a $k$-{\it generic
section},
$k\geq 2$, of an (essential) aspherical arrangement, $\widehat{\cal
A}$. That is, if
there exists $\widehat{\cal A}$ in $\P (V)$, $U\subset V$, with
$M(\widehat{\cal A})$ aspherical, such that ${\cal A}={\widehat{\cal
A}}^U$,
and with the property that $U$ is ${\cal L}_k(\widehat{\cal
A})$-generic, as
in $(1)$ above.

Indeed, Proposition 14 guarantees that in this case we may replace, up
to homotopy,
the inclusion $M({\cal A})\hookrightarrow M(\widehat{\cal A})$ by
a cellular map between minimal CW-complexes, $j: X\to Y$, which
restricts to the
identity on $k$-skeletons. In particular, $Y$ is a $K(\pi, 1)$ and $j$ is
a
classifying map.

Let us recall now from [PS, Def.~2.7] the {\it order of $\pi_1$\/{\rm -}\/connectivity},
$$p(M) := {\rm sup}  ~ \{ q ~|~ b_r(M) = b_r(Y), ~\hbox{ for all } r\leq q
\} ~.\eqno(5)$$
\noindent Set $p = p(M)$. It follows from Proposition 15 that
$p=\infty$ if and only if
$U=V$. Moreover, Propositions 14 and 15 imply that $k \leq p$ and that
we may actually
construct a classifying map $j$ with the property that $j_{|~X^{(p)}}
=~ {\rm id}$
(with the convention $X^{(\infty)}=X$).

We will also need the $\pi$-equivariant
chain complexes of the universal covers, $\widetilde X$ and $\widetilde
Y$,
associated to the above Morse-theoretic {\it minimal} cell structures:
$$ C_*(\widetilde X) := \{ d_q : H_qX\otimes \Z \pi
\longrightarrow H_{q-1}X
\otimes \Z \pi \}_q ~,\eqno(6)$$
 and
$$ C_*(\widetilde Y) := \{ \partial_q : H_qY\otimes \Z \pi
\longrightarrow H_{q-1}Y
\otimes \Z \pi \}_q ~.\eqno(7)$$
 They have the property that $C_{\leq p}(\widetilde X) =
C_{\leq p}(\widetilde Y)$. (Here we adopt the convention of turning
left
$\Z \pi$-modules into right $\Z \pi$-modules, replacing the action of
$x\in \pi$ by that of $x^{-1}$.)

Denote by $\widetilde{j}: C_*(\widetilde X)\to C_*(\widetilde Y)$ the
$\pi$-equivariant chain map induced by the lift of $j$ to universal
covers,
$\widetilde{j} := \{ \widetilde{j}_q :  H_qX\otimes \Z \pi \to
H_qY\otimes \Z \pi \}_q$,
and by $j_* : H_*X \to H_*Y$ the map induced by $j$ on homology,
$$j_* := \{ j_{*q}: H_qX \to H_qY \}_q.$$ Define $\Z \pi$-linear maps,
$\{ D_q \}_q$, by:
$$ D_q :=~\partial_{q+2}+\widetilde{j}_{q+1} :
(H_{q+2}Y\otimes \Z \pi)\oplus (H_{q+1}X\otimes \Z \pi) \longrightarrow
H_{q+1}Y\otimes \Z \pi~.\eqno(8)$$

\nonumproclaim{Theorem 16} Let ${\cal A}$ be an essential
$k$\/{\rm -}\/generic section of an
essential aspherical arrangement $\widehat{\cal A}${\rm ,} with $k\geq 2$.
Set $M=M({\cal A})${\rm ,}
$\pi=\pi_1 (M)$ and $p=p(M)${\rm ,} and denote by $\widetilde M$ the
universal cover
of $M$. Then\/{\rm :}\/

\ritem{(i)} $M$ is aspherical if and only if $p=\infty$.

\ritem{(ii)} If $p<\infty${\rm ,} then the first higher nonzero homotopy group of $M$
is
$\pi_p(M)${\rm ,} with the following finite presentation as a $\Z \pi$\/{\rm -}\/module\/{\rm :}

\vglue-8pt
$$\pi_p(M) = {\rm coker} ~\{ D_p : (H_{p+2}M(\widehat{\cal A})\oplus
H_{p+1}M({\cal A}))
\otimes \Z \pi \longrightarrow H_{p+1}M(\widehat{\cal A})\otimes \Z \pi
\} ~.\eqno(9) $$

\ritem{(iii)} If $2k\geq {\rm rank}({\cal A})${\rm ,} then $\widetilde M$ has the
rational homotopy
type of a wedge of spheres. 

\endproclaim

{\it Proof.} Proposition 15 may be used to infer that both properties
from Part  (i)
of the theorem are equivalent to $U=V$.

The equality $j_{|~X^{(p)}}=~{\rm id}$ implies that $\widetilde M$ is
$(p-1)$-connected; see\break
[PS, Th.~2.10(1)] for details. To deduce the presentation $(9)$ of
$\pi_p(M)$,
we may use the proof of Theorem 2.10(2) from [PS]. There is an
oversight in [PS],
namely the fact that, in general, $\widetilde{j}_{p+1}$ is different
from
$j_{*(p+1)}\otimes {\rm id}$. (Obviously, $\widetilde{j}_q =j_{*q}\otimes
{\rm id}$,
whenever $X^{(q)}$ is a subcomplex of $Y^{(q)}$, and the restriction of
$j$ to  $X^{(q)}$
is the inclusion,  $X^{(q)}\subset Y^{(q)}$.) This is the reason why
$\Delta_p$ from
[PS, (2.3)] has to be replaced by $D_p$ from $(9)$ above. Nevertheless,
this change does not affect the results from [PS, Cor.~2.11], since
$\partial_{p+2}$ is $\varepsilon$-minimal. In particular, $\pi_p(M)\neq
0$, and the proof of
Part (ii) of our Theorem 16 is completed.

 For the last part, let us begin by noting that $\widetilde M$
is a $(k-1)$-connected complex, of dimension $r-1$, $r := {\rm rank}({\cal
A})$.
(The connectivity claim follows from Parts  (i)--(ii), due to the
already remarked fact that
$k\leq p$.) Now the assumption $2k\geq r$
implies firstly that the cohomology algebra $H^*(\widetilde M ; \Q)$
has trivial product structure, and secondly that the rational
homotopy type of $\widetilde M$ is determined by its rational
cohomology algebra
(see [HS, Cor.~5.16]). It is also well-known that the last property
holds
for wedges of spheres; see for example [HS, Lemma 1.6], which finishes
our proof.\eject

\demo{{R}emark {\rm 17}} Denote, as usual, by $M' :=M'({\cal A})$,
the cone of
$M({\cal A})$. Set $\pi' =\pi_1(M')$. The triviality of the Hopf
fibration
readily implies that $M'$ is aspherical if and only if $M$ is
aspherical,
and that $p(M')=p(M)$. Therefore, Theorem 16 (i) also holds for $M'$.
If
$p< \infty$, the Hopf fibration may be invoked once more to see that
$\pi_p(M')=\pi_p(M)$ is the first higher nonzero homotopy group of
$M'$, with
$\Z \pi'$-module structure induced by restriction from $\Z \pi$.
\enddemo

{\it The examples}.
The first explicit computation of nontrivial higher homotopy groups was
made
by Hattori [Hat], for {\it general position} arrangements ${\cal A}$.
Firstly,
he showed that in this case $\pi_1(M({\cal A}))=\Z^n$, where $n= |{\cal
A}|
-1$. Denote by $x_i$ the $1$-cell of the $i$-th $S^1$-factor of the
torus
$(S^1)^n$. One knows that the (acyclic) chain complex $(7)$ coming from
the
canonical (minimal) cell structure of the $n$-torus is of the form
$$
\eqalignno{&\{ \partial_q : \wedge^q (x_1, \dots , x_n) \otimes \Z \Z^n
\longrightarrow
\wedge^{q-1} (x_1, \dots , x_n) \otimes \Z \Z^n \}_q ~,&(10)\cr &&\cr
&\partial_q(x_{i_1}\cdots x_{i_q}) = \sum_{r=1}^q (-1)^{r-1}
x_{i_1}\cdots
\widehat{x_{i_r}}\cdots x_{i_q} \otimes (x_{i_r}^{-1} -1)~.\cr}
$$

Hattori found out that a certain truncation of the above complex
provides a
$\Z \Z^n$-resolution for the first higher nontrivial homotopy group of
$M({\cal A})$.

The second class of examples, studied by Randell in [R1], consists of
{\it iterated generic hyperplane sections}, ${\cal A}=\widehat{\cal
A}^U$,
of essential {\it aspherical} arrangements~$\widehat{\cal A}$. It is
immediate
to see that these are particular cases of (essential) $k$-generic
sections,
with $k= {\rm dim}\ U-1$. Randell gave a formula for the\break $\pi_1$-coinvariants
of the
first higher nontrivial homotopy group of $M'({\cal A})$ (without
determining
the full $\Z \pi_1$-module structure).

The last known examples are the {\it hypersolvable} arrangements,
introduced
in [JP1]. This class contains Hattori's general position class, and
also the
(aspherical) {\it fiber-type} arrangements of Falk and Randell [FR];
see [JP1].
It follows from [JP1, Def.~1.8 and Prop.~1.10] and [JP2, Cor.~3.1] that the (essential) hypersolvable arrangements coincide with
$2$-generic sections
of (essential) fiber-type arrangements. Within the hypersolvable class,
a detailed analysis of the structure of $\pi_2$
as a $\Z \pi_1$-module was carried out in [PS, \S 5].
See also [PS, Th.~4.12(3)]
for a combinatorial formula describing  the $\pi_1$-
coinvariants of the first higher nontrivial homotopy group of the
complement of
hypersolvable arrangements.

Beyond the first higher nontrivial homotopy group, it seems appropriate
to point
out here that Theorem 16 (iii) indicates the presence of a very rich
higher
homotopy structure. Indeed, there are many rank $3$ hypersolvable
arrangements
${\cal A}$ for which $\pi_2M({\cal A})$ is an infinitely generated free
abelian
group; see [PS, Rem. 5.5 and Th.~5.7]. It follows that
$\widetilde{M({\cal A})}$ is rationally an infinite wedge of
$2$-spheres.
Following [Q2], the rational homotopy Lie algebra (under Whitehead
product)
$\pi_{>1}M({\cal A})\otimes \Q$ is a free graded Lie algebra on
infinitely
many generators [Hil].

\vglue-4pt
\section{Some higher homotopy groups of arrangements}
\vglue-6pt

In this section, we are going to apply Theorem 16 to
iterated generic hyperplane sections of aspherical arrangements,
in Theorem 18 below. When specialized to the hypersolvable class,
this result will enable us to prove that the associated graded module
of the
first higher nonvanishing homotopy group of the complement is
combinatorially determined.

\nonumproclaim{Theorem 18} Let $\widehat{\cal A}$ be an essential
aspherical arrangement in $\P ^{m-1}${\rm ,} and let ${\cal A}={\widehat{\cal
A}}^U$
be an iterated generic hyperplane section. {\rm (}\/Here{\rm ,}   a
hyperplane
$H$ is generic if it is ${\cal L}_k$\/{\rm -}\/generic{\rm ,} with $k= {\rm dim}\ H -1${\rm ;} see
\S $5(1)${\rm .)}
Assume ${\rm dim}\ U \geq 3$. Set $M=M({\cal A})${\rm ,} $M'=M'({\cal A})$, $\pi=
\pi_1(M)${\rm ,}
$\pi' =\pi_1(M')${\rm ,} and $p= p(M)= p(M')$. Then{\rm ,}  

\ritem{(i)} $p= {\rm dim}\ U-1$.

\ritem{(ii)} The first higher nontrivial homotopy group of $M$ is $\pi_p(M)${\rm ,}
with
finite{\rm ,} free $\Z \pi$\/{\rm -}\/resolution 

\vglue-12pt
$$\eqalignno{ 0&\to H_{m-1}M(\widehat{\cal A})\otimes \Z \pi \to \cdots&(11)\cr
&\cdots \to
H_{p+2}M(\widehat{\cal A})\otimes \Z \pi \to H_{p+1}M(\widehat{\cal A})
\otimes
\Z \pi \to \pi_p(M) \to 0~,\cr}$$

\item{} where $\{ \partial_q \!: \!H_q M(\widehat{\cal A})\otimes \Z
\pi \to
H_{q-1}M(\widehat{\cal A})\otimes \Z \pi \}_q$ is the
acyclic $\pi$\/{\rm -}\/equivariant chain complex of the universal cover{\rm ,}
associated
to a Morse\/{\rm -}\/theoretic  minimal cell decomposition of $M(\widehat{\cal
A})${\rm ,} as in Section {\rm 5(7)}.
In addition{\rm ,} the resolution $(11)$ is   minimal\/{\rm ;}   more precisely{\rm ,}
$\partial_q \otimes_{\Z \pi} \Z = 0${\rm ,} for all $q${\rm ,} where $\varepsilon :\Z
\pi \to \Z$
is the augmentation of the group ring.

\ritem{(iii)} $\pi_p(M)$ is a projective $\Z \pi$\/{\rm -}\/module $\Longleftrightarrow$
it is $\Z \pi$\/{\rm -}\/free $\Longleftrightarrow {\rm dim}\ U= m-1$.

\ritem{(iv)} The first higher nontrivial homotopy group of $M'$ is
$\pi_p(M')= \pi_p(M)${\rm ,} with $\Z \pi'$\/{\rm -}\/module structure induced from $\Z
\pi${\rm ,}
by restriction of scalars{\rm ,} $\Z \pi' \to \Z \pi${\rm ,} via the projection map
of
the Hopf fibration.

\ritem{(v)} The universal cover of $M$ has the homotopy type of a wedge of
$p$\/{\rm -}\/spheres.
The  rational homotopy Lie algebra{\rm ,} $L_* = \oplus_{q\geq 1}~~L_q${\rm ,}
where
$$L_q := \pi_{q+1}(M)\otimes \Q = \pi_{q+1}(M')\otimes \Q$$ and the Lie
bracket
is induced by the Whitehead product{\rm ,} is isomorphic to the free graded
Lie
algebra generated by $L_{p-1}=\pi_p(M)\otimes \Q$. The\break $\Q\pi$\/{\rm -}\/module
structure
of $L_*$ is given by the graded Lie algebra extension of the
$\pi$\/{\rm -}\/action on
the Lie generators{\rm ,} described by $(11)\otimes \Q${\rm ;} the $\Q\pi'$\/{\rm -}\/module
structure
is obtained by restriction of scalars.

\endproclaim

{\it Proof.} We have already remarked, at the end of the preceding
section,
that ${\cal A}$ is in particular an (essential) $k$-generic section of
$\widehat{\cal A}$, $k= {\rm dim}\ U -1\break \geq 2$, as in Theorem 16. We also
know from
[PS, \S 4.14] that $p= {\rm dim}\ U-1$. Since $H_{p+1}M({\cal A})=0$,
$\,\widetilde{j}_{p+1}=0$. Therefore, the  presentation (9) may be
continued to
the free resolution (11); see (8). The algebraic minimality claim from
Theorem~18(ii) 
follows at once from the topological minimality property of the cell
structure
of $M(\widehat{\cal A})$. Part (iii) is then an immediate
consequence of the minimality of (11), by a standard argument in
homological
algebra (exactly as in [PS, Th.~5.7(1)]). Part (iv) follows
from the triviality of the Hopf fibration.
\vglue6pt
(v) The claim on the homotopy type of the universal cover,
$\widetilde M$,
may be obtained from the remark that $\overline{H}_*(\widetilde{M} ;
\Z)$
is a free abelian group, concentrated in degree $*=p$. This in turn is
a
consequence of the following facts: $\widetilde M$ is
$(p-1)$-connected;
$M$ has the homotopy type of a $p$-dimensional complex. The structure
of
$L_*$ is then provided by Hilton's theorem [Hil]. The $\pi_1(M)$-action
on
$L_*$ is easily seen to be as stated, since one knows that it respects
the Whitehead product [S, p.~419].

\demo{{R}emark {\rm 19}} Hattori's general position framework
corresponds to
the case when $\widehat{\cal A}$ is boolean (that is, when
$\widehat{\cal A}'$
is given by the $m$ coordinate hyperplanes in $\C^m$). See [Hat, \S 2].
In
this case, one recovers, from  (i)--(iv)  above, Hattori's Theorem 3
[Hat].

Randell's formula for the $\pi'$-coinvariants of $\pi_p(M')$ ([R1,$\!$ Th.$\!$~2
and Prop.$\!$~9]) readily follows from (i), (ii) and (iv). He also
raised the question of  $\Z \pi'$-freeness for $\pi_p(M')$, at the end
of
Section 2 [R1]. With the obvious remark that no nonzero $\Z
\pi'$-module,
which is induced by restriction from $\Z \pi$, can be free (since it
has
nontrivial annihilator), (iii) and (iv) above completely clarify
this
question.
\enddemo
 
\centerline{\bf The associated graded chain complex of the universal cover}

\vglue8pt

Let $X$ be a connected CW-complex of finite type, with fundamental
group
$\pi :=\pi_1(X)$, and universal cover $\widetilde X$. The explicit
determination of
the boundary maps, $d_q$, of the equivariant chain complex of
$\widetilde X$, is a
difficult task, in general. Suppose now that $X$ is minimal, with
cohomology ring generated
in degree one. Under these assumptions, we are going to show that,
at the associated graded level, ${\rm gr}^*_I(d_q)$ may be described solely
in terms of the cohomology ring, $H^*X$. When $X$ is (up to homotopy)
an
arrangement complement (either $M({\cal A})$ or $M'({\cal A})$), one
knows [OS] that the
cohomology ring, $H^*X$, is {\it combinatorially determined}, from the
intersection lattice
${\cal L}({\cal A})$. This fact will lead to purely combinatorial
computations of associated graded objects.

Denote by $I\subset \Z \pi$ the augmentation ideal of the group ring,
and by
${\rm gr}^*_I \Z \pi := \oplus_{k\geq 0} (I^k/I^{k+1})$ the associated graded
ring,
with respect to the $I$-adic filtration of $\Z \pi$. Similarly,
we may construct the associated graded module of any (right) $\Z
\pi$-module, $P$,
by ${\rm gr}^*_I P := \oplus_{k\geq 0} (P\cdot I^k/P\cdot I^{k+1})$; it has
an induced
(right) ${\rm gr}^*_I \Z \pi$-module structure.

The minimality of $X$ implies that $d_q(H_qX\otimes \Z \pi)\subset
H_{q-1}X\otimes I$,
for all $q$. Therefore, by passing to the associated graded in $(6)$,
we get a chain complex of free graded  ${\rm gr}^*_I \Z \pi$-modules and
degree zero
${\rm gr}^*_I \Z \pi$-linear maps,
$${\rm gr}^*_I~C_*(\widetilde X) :=~\{ {\rm gr}^*_I(d_q) : H_qX\otimes
{\rm gr}^*_I \Z \pi
\longrightarrow H_{q-1}X\otimes {\rm gr}^*_I \Z \pi \}_q~,\eqno(12) $$
\noindent where we consider on $H_qX\otimes {\rm gr}^*_I \Z \pi$
the tensor product grading, with $H_qX$ concentrated in degree $q$.

Consider now the cup-product map,
$$ H^{q-1}X\otimes H^1X \longrightarrow
H^qX~.\eqno(13)$$
\noindent Dualizing, we get a map, $\delta^R_q : H_qX \to
H_{q-1}X\otimes H_1X$.
Use the standard identifications, $H_1X \equiv \pi_{ab} \equiv I/I^2 =
{\rm gr}^1_I \Z \pi$
(see [HiS, Lemma VI.4.1]), and $H_qX \equiv H_qX\otimes 1$, to extend
it to a
degree zero ${\rm gr}^*_I \Z \pi$-linear map,
$$ \delta^R_q :  H_qX\otimes {\rm gr}^*_I \Z \pi
\longrightarrow H_{q-1}X\otimes {\rm gr}^*_I \Z \pi ~.\eqno(13')$$

\nonumproclaim{Theorem 20} If $X$ is minimal and $H^*X$ is
generated by $H^1X${\rm ,} then 
$$ {\rm gr}^*_I(d_q)~=~(-1)^q \delta^R_q~,~~for~~all~~q~.$$
\endproclaim

{\it Proof.} Set $n := {\rm rank} (\pi_{ab})$. Denote by $Y$ the $n$-torus,
endowed with the canonical minimal CW-structure, giving rise to the
$\pi_{ab}$-equivariant boundary maps, $\partial_q$, described in
$(10)$.
Pick a cellular map, $j : X\to Y$, that induces the abelianization
morphism
on fundamental groups, $j_{\#} : \pi \to \pi_{ab}$. Lift $j$ to a
(cellular,
$j_{\#}$-equivariant) map, $\widetilde{j} : \widetilde X \to \widetilde
Y$,
inducing a $\Z j_{\#}$-linear chain map, $\widetilde{j} :
C_*(\widetilde X) \to C_*(\widetilde Y)$.

Passing to the associated graded, we obtain a morphism of graded rings,
${\rm gr}^*_I (j_{\#}) :  {\rm gr}^*_I \Z \pi \to  {\rm gr}^*_I \Z \pi_{ab}$, and a
degree zero
${\rm gr}^*_I (j_{\#})$-linear chain map,
$$\{ {\rm gr}^*_I (\widetilde{j}_q) :  H_qX\otimes {\rm gr}^*_I \Z \pi \to
H_qY\otimes {\rm gr}^*_I \Z \pi_{ab} \}_q~.$$
\noindent Since ${\rm gr}^0_I (\widetilde{j}_q) \equiv j_{*q}$, by
minimality, we infer that:
$$ {\rm gr}^*_I (\widetilde{j}_q)~=~j_{*q}\otimes  {\rm gr}^*_I
(j_{\#})~,~~for~~all~~q~.\eqno(14)$$

To identify the action of ${\rm gr}^*_I(d_q)$ on the free  ${\rm gr}^*_I \Z
\pi$-module
generators, ${\rm gr}^*_I(d_q) : H_qX \to H_{q-1}X\otimes \pi_{ab}$,
consider  ${\rm gr}^*_I(\partial_q) : H_qY \to H_{q-1}Y\otimes \pi_{ab}$.
Recall that ${\rm gr}^*_I (\widetilde j)$ is a chain map, and use $(14)$
to infer that $(j_{*(q-1)}\otimes j_{*1})\circ {\rm gr}^*_I(d_q) =
{\rm gr}^*_I(\partial_q)\circ j_{*q}$.

At the same time, it is straightforward to use $(10)$ and check that\break
$(-1)^q {\rm gr}^*_I(\partial_q)_{| H_qY}$ equals the dual of the cup-product
map $(13)$,
corresponding to $Y$. By naturality, it follows that
$(j_{*(q-1)}\otimes j_{*1})\circ ((-1)^q \delta^R_q)$ equals
${\rm gr}^*_I(\partial_q)\circ j_{*q}$, as well.

We also know that $H^*X$ is generated by $H^1X$. Since $j_{*1} = {\rm id}$,
this implies that $j$ induces a split injection on homology, and
consequently ${\rm gr}^*_I(d_q) = (-1)^q \delta^R_q$, as asserted.

\vglue12pt

\centerline{\bf Some combinatorially determined situations}

\vglue8pt

In arrangement theory, we may replace ${\rm gr}^*_I \Q \pi$ in Theorem 20
by a combinatorially defined object, in two steps. The first one
uses a general result of Quillen [Q1]. Let $\pi$ be an arbitrary group,
with descending central series $\{ \Gamma_k \}_{k\geq 1}$, and
associated
graded Lie algebra ${\rm gr}^*_{\Gamma} \pi := \oplus_{k\geq 1}
(\Gamma_k/\Gamma_{k+1})$.
One has a map, $\chi_k : \Gamma_k \to I^k$, defined by
$\chi_k (x)=x-1$, for each $k\geq 1$. The maps $\{\chi_k \}$
induce a graded algebra surjection,
$$ \chi : U {\rm gr}^*_{\Gamma} \pi \longrightarrow
{\rm gr}^*_I \Z \pi~,\eqno(15)$$
 which restricts, in degree one, to the standard
identification between
${\rm gr}^1_{\Gamma} \pi = \pi_{ab}$ and ${\rm gr}^1_I \Z \pi$. Moreover, $\chi
\otimes \Q$
is an isomorphism.

The second step involves the {\it holonomy Lie algebra} of $X$, ${\cal
H}_*(X)$.
It is defined as the quotient of the free Lie algebra on $H_1X$, $\L^*
(H_1X)$,
graded by bracket length, by the Lie ideal generated by the image of
the reduced diagonal,
$$ \nabla : H_2 X \longrightarrow H_1X \wedge H_1X \equiv
\L^2 (H_1X)~.\eqno(16)$$
 If $X \simeq M({\cal A})$, one knows [K] that the graded Lie
algebras
$({\rm gr}^*_{\Gamma} \pi)\otimes \Q$ and ${\cal H}_*(X)\otimes \Q$ are
isomorphic,
by an isomorphism which is the identity on degree one pieces.

\nonumproclaim{{C}orollary 21} Let ${\cal A}$ be
an arbitrary projective hyperplane arrangement{\rm ,} with complement
$M({\cal A})$ and
fundamental group $\pi := \pi_1(M({\cal A}))$. Let $X$ be an arbitrary
minimal
complex having the homotopy type of $M({\cal A})$.

\ritem{(i)} The graded algebras  ${\rm gr}^*_I \Q \pi$ and $U {\cal H}_*(M({\cal
A}))\otimes \Q$
are isomorphic{\rm ,} by an isomorphism inducing the identity in degree one.
Consequently{\rm ,}
the chain complex of graded  ${\rm gr}^*_I \Q \pi$\/{\rm -}\/modules ${\rm gr}^*_I
C_*(\widetilde X)\otimes \Q$
is combinatorially determined.

\ritem{(ii)} If the associated central arrangement{\rm ,} ${\cal A}'${\rm ,} is
hypersolvable{\rm ,} then
all the results from Part {\rm (i)} hold with $\Z$\/{\rm -}\/coefficients. 

\endproclaim

{\it Proof.} Part (i) is a direct consequence of Theorem 20, via the
results from
[Q1] and [K] mentioned above. For Part (ii), it will be plainly
enough to establish
the asserted isomorphism between  ${\rm gr}^*_I \Z \pi$ and   $U {\cal
H}_*(M({\cal A}))$.

To this end, we will consider the complement of ${\cal A}'$, denoted,
as usual, by
$M'({\cal A})$, with fundamental group $\pi' := \pi_1(M'({\cal A}))$.
Up to homotopy,
$M'({\cal A})\simeq M({\cal A}) \times S^1$. We are going to exploit
this simple remark,
to obtain the result we need for $M({\cal A})$ from known facts about
$M'({\cal A})$.

Since the group $\pi$ is a retract of the group $\pi'$, the graded Lie
algebra
${\rm gr}^*_{\Gamma} \pi$ embeds into ${\rm gr}^*_{\Gamma} \pi'$. One also knows
that the latter
is a free abelian group; see [JP1, Th.~C]. It follows that $U
{\rm gr}^*_{\Gamma} \pi$ is
a free abelian group as well, by Poincar\'{e}-Birkhoff-Witt, and
consequently the Quillen map
$(15)$ is an isomorphism.

One also knows [MP, Prop.~5.1] that there is a graded Lie algebra
surjection,
inducing the identity in degree one,
$$ \phi_S~:~{\cal H}_*(S) \longrightarrow
{\rm gr}^*_{\Gamma} \pi_1(S)~,\eqno(17)$$
\noindent for any connected space $S$ with finite Betti numbers and
without torsion in
$H_1S$. The result from [K] quoted above implies then that
$\phi_{M({\cal A})}\otimes \Q$
is an isomorphism. To show that actually $\phi_{M({\cal A})}$ is an
isomorphism, and
thus finish the proof of our corollary, it will be enough to check that
${\cal H}_*(M({\cal A}))$
is a free abelian group. To do this, we may argue as before, starting
with the remark that
$M({\cal A})$ is a retract of $M'({\cal A})$. By naturality (see
definition $(16)$), this implies that
${\cal H}_*(M({\cal A}))$ embeds into ${\cal H}_*(M'({\cal A}))$.
Finally,  ${\cal H}_*(M'({\cal A}))$
is a free abelian group, by [JP1, Th.~C] (see also [JP1, (7.6)]),
which completes our proof.

\vglue12pt

To obtain our last main result, Theorem 23, we will use the Koszul
property as key ingredient. So
we start by recalling from [Lof] the definition of Koszulness.

Let $U^*$ be a connected graded algebra of finite type, over a field
$\K$. The algebra $U$ is
{\it Koszul} if ${\rm Tor}^U_{s,t} (\K, \K)=0$ for $s\neq t$, where $s$ is
the homological (resolution) degree,
and $t$ is the internal degree (coming from the grading of $U^*$).

The Koszul property from Part (i) of the technical proposition below
will be used to obtain an exactness property at the associated graded
level, in the (aspherical) cases treated in Part (ii).
This result in turn will provide an important step in the proof of
Theorem 23.

\nonumproclaim{Proposition 22} Let ${\cal A}$ be a projective
arrangement whose cone{\rm ,} ${\cal A}'${\rm ,}
is hypersolvable. Set $M:=M({\cal A})${\rm ,} and $U^*:= U {\cal H}_*(M)$.
Let $X$ be any minimal cell
structure of $M$.

\ritem{(i)} The algebra $U\otimes \K$ is Koszul{\rm ,} for any field $\K$.

\ritem{(ii)} Assume that ${\cal A}'$ is fiber\/{\rm -}\/type. Then the augmentation{\rm ,}
$\varepsilon : U\to \Z${\rm ,} extends to
a free {\rm (}\/right\/{\rm )} $U^*$\/{\rm -}\/resolution{\rm ,} $\varepsilon : {\rm gr}^*_I C_*(\widetilde
X)\to \Z${\rm ,}   where\/{\rm :}\/

\vglue-8pt $$ {\rm gr}^*_I C_*(\widetilde X)~=~\{ (-1)^q \delta^R_q :
H_qM\otimes U^* \longrightarrow H_{q-1}M\otimes U^* \}_q~,\eqno(18)$$

\item{} and the restriction{\rm ,}
$\delta^R_q : H_qM\otimes 1 \equiv H_qM \to H_{q-1}M\otimes U^1 \equiv
H_{q-1}M\otimes H_1M${\rm ,}
is dual to the cup\/{\rm -}\/product{\rm ,} $H^{q-1}M\otimes H^1M \to H^qM$.

\endproclaim

{\it Proof.} (i) This will follow from the fact that $M$ is a retract
of $M':=M'({\cal A})\simeq M\times S^1$,
as in the proof of Corollary 21 (ii). By naturality, the graded
algebra
$U\otimes \K =  U {\cal H}_*(M)\otimes \K$ is a retract of
$U'\otimes \K :=  U {\cal H}_*(M')\otimes \K$, and therefore
${\rm Tor}^{U\otimes \K}_{*,*}(\K, \K)$ is a bigraded $\K$-
subspace of ${\rm Tor}^{U'\otimes \K}_{*,*}(\K, \K)$. We deduce from [JP1,
Th.~E] that $U'\otimes \K$
is a Koszul algebra, via Koszul duality (see [JP1, 7(2)]). The
Koszulness of $U\otimes \K$ follows.

(ii) Everything is a direct consequence of Theorem 20 and Corollary
21 (ii), except
the fact that the augmented chain complex $(18)$ is acyclic. It is
enough to check the acyclicity property
for arbitrary field coefficients $\K$.

As a preliminary step, we will prove first that the algebra
$H^*M\otimes \K$ is Koszul. As before,
the topological retraction property implies the existence of a split
graded algebra epimorphism,
$H^*M'\otimes \K \to H^*M\otimes \K$, inducing a surjection on
${\rm Tor}^{(\cdot)}_{*,*}(\K, \K)$.
Since $H^*M'\otimes \K$ is a Koszul algebra, for the fiber-type class
(by the main result of [SY]),
we infer that $H^*M\otimes \K$ is also Koszul.

At this point, we will need the {\it Koszul duality}. Let $A^* =
\T^*(V)/ {\rm ideal}~(R)$
be a {\it quadratic algebra}, where the tensor algebra $\T^*(V)$ is
graded by tensor degree,
and $R$ is a $\K$-vector subspace of $V^{\otimes 2}$. Denote by $\#
(\cdot)$ $\K$-duals
of vector spaces. The {\it Koszul dual} of $A^*$ is the quadratic
algebra
$A^{!*} := \T^*(\# V)/ {\rm ideal}~(R^{\perp})$, where $R^{\perp}\subset \#
V^{\otimes 2}$
is the annihilator of $R$.

 From the preceding definition and definition $(16)$, and by
the fact that
$H^*(M, \K)=H^*M\otimes \K$ is a (strictly) graded-commutative algebra,
it is not difficult to verify that the graded algebras $(U\otimes
\K)^{!*}$ and
$\wedge^*_{\K}(H^1M)/{\rm ideal}~({\rm ker}~\mu)$ are isomorphic. Here
$\wedge^*_{\K}(H^1M)$
denotes the exterior $\K$-algebra on $H^1M\otimes \K$, graded by
exterior degree,
and $\mu : \wedge ^2_{\K}(H^1M)\to H^2M\otimes \K$ is the
cup-multiplication.

The vanishing property of the ${\rm Tor}$-pieces for $s=1$ and $s=2$,
guaranteed by the Koszulness of $H^*M\otimes \K$,
translates to the fact that the algebra $H^*M\otimes \K$ is generated
by $H^1M\otimes \K$,
and the defining relations are of degree $2$. In other words,
$H^*M\otimes \K = (U\otimes \K)^{!*}$,
as graded algebras.

The Koszul property of $U\otimes \K$ from Part (i) implies, via the
above identification,
that the augmentation, $\varepsilon : U\otimes \K \to \K$, extends to a
(left) free
$U^*\otimes \K$-resolution, $\varepsilon : C^L_* \to \K$, where
$$ C^L_* := \{ - \delta^L_q : (U^*\otimes \K)\otimes H_q(M,
\K)\longrightarrow
(U^*\otimes \K)\otimes H_{q-1}(M, \K) \}_q~,\eqno(19)$$
 and the restriction, $\delta^L_q : H_q(M, \K)\to H_1(M,
\K)\otimes H_{q-1}(M, \K)$,
is dual to the cup-product map,
$$ H^1(M, \K)\otimes H^{q-1}(M, \K)
\longrightarrow H^q(M, \K)~.\eqno(20)$$
  See [Lof, pp.\ 305--306]  for a proof of the fact that the
above
{\it resolution} property actually characterizes Koszul algebras,
within
the quadratic class.

By construction, the graded algebra $U^*\!\otimes \K$ is the quotient of
$\T^*(H_1(M, \K))$ by the ideal generated by ${\rm im}~(\nabla \otimes
\K)\subset \T^2$.
The involutive graded algebra anti-automorphism, $\alpha$, of $\T^*$,
which sends
$x_1\otimes \cdots \otimes x_k \in \T^k$ to $x_k\otimes \cdots \otimes
x_1$,
induces a graded algebra anti-automorphism, $\alpha : U^*\!\otimes \K \to
U^*\otimes \K$.
Denote by $c: V\otimes W \to W \otimes V$ the $\K$-linear isomorphism
which
interchanges the factors of a tensor product of two $\K$-vector spaces.
Finally,
consider the $\K$-isomorphisms
$$\{ ({\rm id} \otimes \alpha)\circ c~:~(U^*\otimes \K)\otimes H_q(M,
\K)\longrightarrow
H_q(M, \K)\otimes (U^*\otimes \K) \}_q~.$$
\noindent It is easy to check that they define a $\K$-linear
isomorphism
between the augmented chain complexes $(18)\otimes \K$ and $(19)$.
Therefore,
$\varepsilon : {\rm gr}^*_I C_*(\widetilde X)\otimes \K \to \K$ is acyclic, and
this finishes the proof of Proposition 22.

Let ${\cal A}$ be an essential projective arrangement of rank $r\geq
3$, with
{\it hypersolvable} cone, ${\cal A}'$. Set $M:=M({\cal A})$, and
$p:=p(M)$.
Theorem 18 applies to $M$ (that is, ${\cal A}$ is an iterated generic
hyperplane section
of an essential aspherical arrangement $\widehat{\cal A}$) precisely
when $p=r-1$.
In one direction, the above claim follows from Theorem 18 (i). In
general, one knows [JP2]
that ${\cal A}$ is a $2$-generic section of an essential fiber-type
arrangement,
$\widehat{\cal A}$. Assume now that $p=r-1$. We may then infer from
Proposition 15,
via our discussion preceding Theorem~16, on the $k$-generic framework,
that actually ${\cal A}$
is a $p$-generic section of~$\widehat{\cal A}$, i.e., an iterated
generic hyperplane section of
$\widehat{\cal A}$.

The set of (essential) hypersolvable arrangements with $p=r-1$ is
combinatorially defined.
Indeed, both the definition of the hypersolvable class and the
computation of $p$ within this class
are purely combinatorial; see [JP1, Def.~1.8] and [PS, Cor.~4.10(1)]
respectively.
(Note also that the condition $r\geq 3$ is automatically fulfilled,
since $p\geq 2$; see
[PS, Cor.~4.10(2)].)

The hypersolvable examples with $p=r-1$ abound. Let us consider for
instance the rank $3$
hypersolvable case. Here, $p=2$ if and only if $\ell >3$ (by [PS, Th.~5.4(1)]), where
$\ell =\ell ({\cal A}')$ is the length of ${\cal A}'$, introduced in
[JP1]. At the same time,
rank $3$ hypersolvable examples of arbitrarily large length are easily
constructed; see
[JP1, \S 1]. One may also obtain hypersolvable examples with
arbitrarily large
rank $r$, and $p=r-1$, by taking generic hyperplane sections of
fiber-type arrangements
of large rank; see [PS, p.~90].

To formulate our next result, we need to recall one more fact, from\break
[JP1, Prop.~3.2(i)]:
one may associate to any hypersolvable arrangement~${\cal A}'$ a
combinatorially determined
collection of positive natural numbers,\break $\{ 1=d_1, d_2,\dots ,d_{\ell}
\}$, which are
called the {\it exponents} of ${\cal A}'$.

\nonumproclaim{Theorem 23} Let ${\cal A}$ be an essential
projective arrangement with
hypersolvable cone{\rm ,} ${\cal A}'$. Set $M:=M({\cal A})${\rm ,} $\pi:=\pi_1(M)${\rm ,}
$p:=p(M)${\rm ,} $r:={\rm rank} ({\cal A})$.
Assume that $p=r-1${\rm ,} and realize ${\cal A}$ as an iterated generic
hyperplane section
of an essential fiber\/{\rm -}\/type arrangement{\rm ,} $\widehat{\cal A}$.
\eject
 {\rm (i)} There is an isomorphism of graded rings{\rm ,} ${\rm gr}^*_I \Z \pi \equiv
U^*${\rm ,}
where $U^* := U {\cal H}_*(M(\widehat{\cal A}))$ and ${\cal H}_*$
denotes the holonomy Lie algebra
defined by  {\rm (16),} and a compatible isomorphism of graded modules{\rm ,}
$$ {\rm gr}^*_I \pi_p(M) \equiv \Sigma^{-p-1}~{\rm coker}~\{ \delta^R_{p+2}
: H_{p+2}M(\widehat{\cal A})
\otimes U^* \longrightarrow H_{p+1}M(\widehat{\cal A})\otimes U^*
\}~,\eqno(21)$$
  where $\delta^R_{p+2}$ is constructed as in Proposition
 {\rm 22(ii),} and $\Sigma^{-p-1}$
denotes the\break $(p+1)$\/{\rm -}\/desus\/{\rm -}\/pension operation on graded modules. In
particular{\rm ,} the graded module
${\rm gr}^*_I \pi_p(M)$ is combinatorially determined. More precisely{\rm ,} it
depends only on
the rank of ${\cal A}${\rm ,} and on ${\cal L}_2({\cal A})${\rm ,}
the elements of rank at most $2$ from ${\cal L}({\cal A})$.

\vglue4pt {\rm (ii)} All graded pieces{\rm ,} ${\rm gr}^i_I \pi_p(M)${\rm ,} are finitely\/{\rm -}\/generated free
abelian groups. The
Hilbert series{\rm ,} ${\rm gr}^*_I \pi_p(M)(t) := \sum_{i\geq 0}~{\rm rank}~{\rm gr}^i_I
\pi_p(M)\cdot t^i${\rm ,} is equal to
$$ (-1/t)^{p+1}~\{ 1- [\sum_{j=0}^p~(-1)^j \beta_j \cdot
t^j]/
\prod_{i=2}^{\ell}~(1-d_i t) \}~,\eqno(22)$$
 where  $\{ 1=d_1, d_2,\dots ,d_{\ell} \}$ are the exponents
of ${\cal A}'${\rm ,} and
the coefficients\break $\{ \beta_j = b_j(M(\widehat{\cal A})) \}$ are given
by
$$ (1+d_2 t) \cdots (1+d_{\ell} t)~=~\sum_{j=0}^{\ell
-1}~\beta_j \cdot t^j~.\eqno(23)$$
 In particular{\rm ,}  ${\rm gr}^*_I \pi_p(M)(t)$ depends only on the rank
and the exponents of ${\cal A}'$.
\endproclaim

{\it Proof.} (i) The $\Z \pi$-module $\pi_p(M)$ is isomorphic to
$${\rm coker}~\{ \partial_{p+2} : H_{p+2}M(\widehat{\cal A})
\otimes \Z \pi  \to H_{p+1}M(\widehat{\cal A})\otimes \Z \pi \}~,$$
by Theorem 18(ii). The following notation will be useful. Denote by
$M_i$ the\break $\Z \pi$-modules
$H_{p+i}M(\widehat{\cal A})\otimes \Z \pi$, for $1\leq i \leq 3$, and
by $\{ F_k M_i \}_{k\geq 0}$
the corresponding $I$-adic filtrations. Set $\partial :=
\partial_{p+2}$, and $\partial' := \partial_{p+3}$.

Next, we shall isolate three facts which are key to our proof.
Firstly,\break $\partial \partial' =0$.
Secondly, $\partial (F_*M_2)\subset F_{*+1}M_1$, and $\partial'
(F_*M_3)\subset F_{*+1}M_2$, by minimality.
Finally,
$$ \eqalignno{&\hskip-.5in{\rm im}~\{ {\rm gr}^*(\partial') : {\rm gr}^{*-1}(M_3) \longrightarrow {\rm gr}^*(M_2)
\}&(24)\cr
&\quad =~
 {\rm ker}~\{ {\rm gr}^{*+1}(\partial) : {\rm gr}^*(M_2) \longrightarrow {\rm gr}^{*+1}(M_1)
\}~.\cr}$$
 This crucial property follows from Proposition 22(ii).

The canonical surjections, $$F_*M_1 / F_{*+1}M_1 \to [F_*M_1 +
~{\rm im}~(\partial)] / [F_{*+1}M_1 +~{\rm im}~(\partial)] =
{\rm gr}^*_I \pi_p(M)~,$$\eject\noindent  factor to give a surjective morphism of graded
${\rm gr}^*_I \Z \pi$-modules,
$$ \psi : ~co {\rm ker}~\{ {\rm gr}^*(\partial) \}
\longrightarrow  {\rm gr}^*_I \pi_p(M)~.\eqno(25)$$
\noindent We will resort to the above mentioned three basic facts and
prove that $\psi$ is injective,
in each degree $*=k$.

This means that, given any $v^1$ such that $v^1 \equiv_k 0$ and $v^1
\equiv_{k+1} \partial u^2$, we
have to find $w^2$ such that $w^2 \equiv_{k-1} 0$ and $v^1 \equiv_{k+1}
\partial w^2$.
Here the superscripts $j$ indicate that the corresponding elements
belong to $M_j$, and the notation $\equiv_q$ means equality modulo the
appropriate
filtration term $F_q$. We will argue by induction and suppose that we
have found
$w^2_s$ such that $w^2_s \equiv_s 0$ and $v^1 \equiv_{k+1} \partial
w^2_s$, for $s < k-1$.
(For $s=0$, we may take $w^2_0=u^2$.) The equalities $\partial w^2_s
\equiv_{s+2} v^1 \equiv_{s+2} 0$
imply the existence of $z^3$ such that $z^3 \equiv_{s-1} 0$ and
$w^2_s \equiv_{s+1} \partial' z^3$,
due to the exactness property $(24)$. Setting $w^2_{s+1} = w^2_s -
\partial' z^3$, it is immediate to check that
$w^2_{s+1} \equiv_{s+1} 0$ and $v^1 \equiv_{k+1} \partial w^2_{s+1}$,
as needed for the induction step.
Finally, we may take $w^2 = w^2_{k-1}$, to complete the proof of the
fact that $\psi$ is a
${\rm gr}^*_I \Z \pi$-isomorphism.

We may now finish the proof of Theorem 23(i) as follows. The asserted
identification of graded rings,
${\rm gr}^*_I \Z \pi \equiv U^*$, is provided by Corollary 21(ii). The
compatible isomorphism of graded modules
 (21)  becomes a consequence of Proposition 22(ii). Therefore, the
graded module  ${\rm gr}^*_I \pi_p(M)$ is
determined by $p$ and the cohomology ring of $M(\widehat{\cal A})$.

To see that $H^*M(\widehat{\cal A})$ is in turn determined by  ${\cal
L}_2({\cal A})$, we may argue as follows.
By Proposition 14, the rings $H^*M(\widehat{\cal A})$ and $H^*M({\cal
A})$ are isomorphic, up to degree $2$.
Since ${\cal L}_2({\cal A})$ determines $H^{\leq 2}M({\cal A})$, by
[OS], we infer that the $\Z$-algebra
$$\wedge^*(H^1M(\widehat{\cal A}))/ {\rm ideal}~( {\rm ker}~\{ \mu :
\wedge^2(H^1M(\widehat{\cal A})) \to H^2M(\widehat{\cal A})\} )$$
  depends only on ${\cal L}_2({\cal A})$. The above algebra
naturally surjects onto $H^*M(\widehat{\cal A})$.
It follows from the proof of Proposition 22(ii) that this surjection
induces an isomorphism,
with arbitrary field coefficients $\K$, due to the Koszul property of
$H^*M(\widehat{\cal A})\otimes \K$.
This implies that the isomorphism actually holds over $\Z$, and we are
done.

(ii) By Part (i) and the resolution property from Proposition
22 (ii),
we deduce that ${\rm gr}^*_I \pi_p(M)$ is isomorphic to
$$\Sigma^{-p-1}~{\rm im}~(\delta^R_{p+1})\subset
\Sigma^{-p-1}(H_pM(\widehat{\cal A})\otimes U^*)~.$$ The first assertion
from Theorem 23(ii)
follows then from the fact that the graded pieces of $U^* = U {\cal
H}_*(M(\widehat{\cal A}))$
are finitely-generated free abelian groups; see the proof of Corollary
21(ii).

To obtain the $I$-adic filtration formula $(22)$, we will use
$\Q$-coefficients and the resolution from
Proposition 22(ii). Denote by $\{B^*_q \}_{q\in \Z}$ and $\{ K^*_q
\}_{q\in \Z}$ respectively,
the $q$-boundaries and the $q$-chains of the acyclic chain complex
$\varepsilon : {\rm gr}^*_I C_*(\widetilde X)\otimes \Q \to \Q$. Use the exact
sequences
$0 \to B^*_q \to K^*_q \to B^*_{q-1} \to 0$ and straightforward
induction to obtain the
following formula for the Hilbert series of $B^*_p$:
$$ (-1)^{p+1} B_p(t)~=~1 + \sum_{j=0}^p (-1)^{j+1}
b_j(M(\widehat{\cal A})) t^j \cdot U\otimes \Q (t)~.\eqno(26)$$
  On the other hand, it follows from $(21)$ that
$$t^{p+1}  {\rm gr}^*_I \pi_p(M)(t)~=~B_p(t)~.\eqno(27) $$
 The resolution property also implies that
$$ H^*(M(\widehat{\cal A}), \Q)(-t) \cdot U\otimes \Q
(t)~=~1~,\eqno(28)$$
  via a standard Euler characteristic argument.

Finally, one also knows that $\widehat{\cal A}'$ and ${\cal A}'$ have
the same exponents; see
[JP2, \S 2]. We infer from [FR] that $H^*(M'(\widehat{\cal A}), \Q)(t)
= \prod_{i=1}^{\ell} (1+ d_i t)$,
and consequently
$$ H^*(M(\widehat{\cal A}), \Q)(t) ~=~
\prod_{i=2}^{\ell} (1+ d_i t)~.\eqno(29)$$

Equations (26)--(29) together establish formula $(22)$, thus finishing
the proof of Theorem 23.

\vglue12pt {\it Acknowledgements}.
The authors thank
Paltin Ionescu, who gave a proof of Corollary 4, Part (1), with
completely
different methods, for $n=2$, which challenged us to look for a general
proof.
He also noticed the relation of our results to Dolgachev's paper [Do];
in particular we owe Corollary 4, Part  (2)(b), to him.
We also thank Pierrette Cassou-Nogu\`es and Alex Suciu for
drawing our attention to Richard Randell's preprint [R2].
The first author thanks Jim Damon and Dirk Siersma for very
useful
discussions on a preliminary version of this paper and for pointing the
relations with the papers [Da1,2,3], [JL], [Si] and [V].

\AuthorRefNames [DPW]
\references

 [B]  \name{E.\ Brieskorn}, Sur les groupes de tresses, in {\it S\'eminaire
Bourbaki\/} 1971/72, {\it Lecture Notes in Math\/}.\ {\bf 317}, 21--44,
Springer-Verlag, New York (1973).  

[CD] \name{Pi.\ Cassou-Nogu\`es} and \name{A.\  Dimca},  Topology of complex
polynomials via polar curves, {\it Kodai Math.\ J\/}.\ {\bf 22} (1999),
131--139.

[C] \name{D.\ C.\  Cohen}, Morse inequalities for arrangements, {\it
Adv.\ Math\/}.\ {\bf 134} (1998), 43--45.

[CS] \name{D.\ C.\  Cohen} and \name{A.\  Suciu}, Homology of iterated semidirect
products
of free groups, {\it J.\ Pure Appl.\ Alg\/}.\ {\bf 126}  (1998), 87--120.

[Da1] \name{J.\ Damon},  {\it Higher Multiplicities and Almost Free Divisors
and Complete Intersections\/}, {\it Mem.\  Amer.\ Math.\ Soc\/}.\ {\bf
123}, Providence, RI (1996).

[Da2] \bibline,  Critical points of affine multiforms on the
complements of arrangements, in {\it Singularity Theory\/} (J.\ W.\
Bruce and D.\ Mond, ed.), {\it London Math.\ Soc.\ Lecture Notes\/} {\bf
263}, 25--53, Cambridge, Univ.\ Press, Cambridge (1999).

[Da3] \bibline, On the number of bounding cycles for nonlinear
arrangements, {\it Arrangements\/} (Tokyo, 1998), {\it Adv.\ Stud.\  Pure
Math\/}.\ {\bf 27} (2000), 51--72.

[DaM] \name{J.\ Damon} and \name{D.\ Mond},  $A$-codimension and the vanishing
topology
of discriminants, {\it Invent.\ Math\/}.\ {\bf 106}  (1991), 217--242.

[De] \name{P.\ Deligne}, Les immeubles des groupes de tresses
g\'en\'eralis\'es, {\it Invent.\ Math\/}.\ {\bf 17}  (1972), 273--302.

[D1] \name{A.\ Dimca}, {\it Singularities and Topology of
Hypersurfaces\/},
{\it Universitext\/}, Springer-\break Verlag, New York, 1992.

[D2] \bibline, Arrangements, Milnor fibers and polar curves,
preprint; math.AG/0011073, version 2.

[D3] \bibline,  On polar Cremona transformations, {\it An.\ St.\
Univ.\ Ovidius Constanta\/} {\bf 9} (2001), 47--54.

 [DL] \name{A.\ Dimca} and \name{G.\ I.\ Lehrer}, Purity and equivariant weight
polynomials, in {\it Algebraic Groups and Lie Groups\/} (G.\ I.\ Lehrer,
ed.), Cambridge Univ. Press, Cambridge (1997), 161--181.

[DP] \name{A.\  Dimca} and \name{L.\  P\u aunescu}, On the connectivity of
complex affine hypersurfaces.\ II, {\it Topology\/} {\bf 39}  (2000),
1035--1043.

[Do] \name{I.\ V.\  Dolgachev}, Polar Cremona transformations, {\it Michigan
Math.\ J\/}.\ {\bf 48} (2000) (volume dedicated to W.\ Fulton), 191--202.

[dPW] \name{A.\ A.\  du Plessis} and \name{C.\ T.\ C.\  Wall}, Discriminants, vector
fields and singular hypersurfaces, in {\it New Developments in Singularity
Theory\/} (D.\ Siersma, et al.\ eds.),  Kluwer Academic Publ., Dordrecht
(2001), 351--377.

[FR] \name{M.\ Falk} and \name{R.\ Randell}, The lower central series of a
fiber-type
arrangement, {\it Invent.\ Math\/}.\ {\bf 82} (1985), 77--88.

[GM] \name{M.\ Goresky} and \name{R.\  MacPherson}, {\it Stratified Morse
Theory\/},
{\it Ergeb.\ Math.\ Grenzgeb\/}.\ {\bf 14}, Springer-Verlag, New York,  1988.

[HS] \name{S.\ Halperin} and \name{J.\ Stasheff}, Obstructions to homotopy
equivalences,
{\it Adv.\  Math\/}.\ {\bf 32}  (1979), 233--279.

[H] \name{H.\ A.\  Hamm}, Lefschetz theorems for singular varieties,
{\it Proc.\ Sympos.\ Pure Math\/}.\ {\bf 40} (1983),
547--557.

[HL] \name{H.\ A.\  Hamm} and \name{D.\ T.\  L\^ e},  Un th\' eor\` eme de Zariski du
type de
Lefschetz, {\it Ann.\ Sci.\ \' Ecole Norm.\ Sup\/}.\ {\bf 6}   (1973),
317--355.

[Hat] \name{A.\ Hattori}, Topology of ${\scriptstyle\C}^n$ minus a finite number of
affine hyperplanes in general position, {\it J.\ Fac.\ Sci.\ Univ.\
Tokyo\/} {\bf 22} 
(1975), 205--219.

[Hil] \name{P.\ J.\  Hilton}, On the homotopy groups of the union of
spheres,
{\it J.\ London Math.\ Soc\/}.\ {\bf 30}  (1955), 154--172.

[HiS] \name{P.\ J.\  Hilton} and \name{U.\  Stammbach},  {\it A Course in Homological
Algebra\/}, {\it  Grad.\ Texts in Math\/}.\ {\bf 4},
Springer-Verlag, New York, 1971.

[Hi] \name{H.\ Hironaka}, Stratifications and flatness, in {\it Real and
Complex
Singularities\/}, 199--265, Sijthoff and Noordhoff, Alphen aan den Rijn (1977).

[JL] \name{Ph.\ Jacobs} and \name{A.\ Laeremans}, Complex arrangements whose
complement has
Euler characteristic $\pm 1$, {\it Indag.\ Math\/}.\ {\bf 8}   (1997),
399--404.

[JP1] \name{M.\ Jambu} and \name{S.\  Papadima},  A generalization of fiber-type
arrangements
and a new deformation method, {\it Topology\/} {\bf 37}  (1998), 1135--1164.

[JP2] \bibline, Deformations of hypersolvable
arrangements,
{\it Topology Appl\/}.\ {\bf 118} (2002), 103--111.

[K] \name{T.\ Kohno}, On the holonomy Lie algebra and the nilpotent
completion of
the fundamental group of the complement of hypersurfaces, {\it Nagoya
Math\/}.\ {\it J\/}.\ 
{\bf 92}  (1983), 21--37.

[L] \name{K.\ Lamotke},  The topology of complex projective varieties
after S. Lefschetz, {\it Topology\/} {\bf 20}  (1981), 15--51.

 [Le1] \name{D.\ T.\ L\^e},  Calcul du nombre de cycles \'evanouissants
d'une hypersurface complexe, {\it Ann.\ Inst.\ Fourier\/}, {\it
Grenoble\/} {\bf 23} (1973),
261--270.

 [Le2] \bibline,   Le concept de singularit\'e isol\'ee de
fonction analytique, in {\it Complex Analytic Singularities\/},
{\it Adv.\ Studies  Pure Math\/}.\ {\bf 8}, 215--227,
North-Holland, Amsterdam (1986).

[LT] \name{D.\ T.\ L\^e} and \name{B.\ Teissier},  Vari\'et\'es polaires locales et
classes de Chern des vari\'et\'es singuli\`eres, {\it Ann.\ of
Math\/}.\ {\bf 114\/} (1981),
457--491.

[Li] \name{A.\ Libgober},  Homotopy groups of the complements to singular
hypersurfaces, II, {\it Ann.\ of  Math\/}.\ {\bf 139} (1994), 117--144.

[Lof] \name{C.\ L\" ofwall}, On the subalgebra generated by
one-dimensional elements in
the Yoneda Ext-algebra, in {\it Algebra\/}, {\it Algebraic Topology and their
Interactions\/}, {\it Lecture Notes in Math\/}.\ {\bf 1883} (1986), 291--338,
Springer-Verlag, New York.

[Lo] \name{E.\ Looijenga},  {\it Isolated Singular Points on Complete
Intersections\/}, {\it London Math.\ Soc.\ Lecture Note Ser\/}.\ {\bf 77}, 
Cambridge, Univ.\ Press, Cambridge, 1984.

[MP] \name{M.\ Markl} and \name{S.\  Papadima}, Homotopy Lie algebras and
fundamental groups via deformation theory,
{\it Ann.\ Inst.\ Fourier\/} {\bf 42}  (1992), 905--935.

 [M] \name{J.\ Milnor}, {\it Singular Points of Complex Hypersurfaces\/},
{\it Ann.\ of Math.\ Studies\/} {\bf 61}, Princeton Univ.\ Press,
Princeton, NJ, 1968.

[MO] \name{J.\ Milnor} and \name{P.\  Orlik}, Isolated singularities defined by
weighted homogeneous polynomials, {\it Topology\/} {\bf 9}  (1970), 385--393.

[N1] \name{A.\ N\'emethi}, Th\'eorie de Lefschetz pour les vari\'et\'es
alg\'ebriques affines, {\it C.\ R.\ Acad.\ Sci.\ Paris S\'er.\ I
Math\/}.\ {\bf 303}  
(1986), 567--570.

[N2] \bibline, Lefschetz theory for complex affine
varieties, {\it Rev.\ Roumaine Math.\ Pures  Appl\/}.\ {\bf 33}  (1988), 
233--250.

[OS] \name{P.\ Orlik} and \name{L.\  Solomon}, Combinatorics and topology of
complements
of hyperplanes, {\it Invent.\ Math\/}.\ {\bf 56}  (1980), 167--189.

[OT1] \name{P.\ Orlik} and \name{H.\  Terao}, {\it Arrangements of
Hyperplanes\/},
{\it Grundlehren Math.\ Wiss\/}.\ {\bf 300}, Springer-Verlag, New York, 1992.

[OT2] \bibline, Arrangements and Milnor fibers,
{\it Math.\ Ann\/}.\ {\bf 301}  (1995), 211--235.

[OT3] \bibline, The number of critical points of a
product of powers of linear functions,
{\it Invent.\ Math\/}.\ {\bf 120}  (1995), 1--14.

[PS] \name{S.\ Papadima} and \name{A.\  Suciu}, Higher homotopy groups of
complements of complex hyperplane arrangements, {\it Adv.\  Math\/}.\
{\bf 165} 
(2002), 71--100.

[Q1] \name{D.\ Quillen}, On the associated graded ring of a group ring,
{\it J.\ Algebra\/} {\bf 10}  (1968), 411--418.

[Q2] \bibline, Rational homotopy theory, {\it Ann.\ of
Math\/}.\ {\bf 90}  (1969),
205--295.

[R1] \name{R.\ Randell}, Homotopy and group cohomology of arrangements,
{\it Topology Appl\/}.\ {\bf 78}  (1997), 201--213.

[R2] \name{R.\ Randell},  Morse theory, Milnor fibers and hyperplane
arrangements, preprint;\break   math.AT 0011101, version 2.

[R3] \bibline, Morse theory, Milnor fibers and minimality of
hyperplane arrangements, {\it Proc.\ Amer.\ Math.\ Soc\/}.\ {\bf 130}  
(2002), 2737--2743.

[Ry] \name{G.\ Rybnikov}, On the fundamental group of the complement of
a complex hyperplane arrangement, DIMACS Tech.\ Report 94-13 (1994),
33-50; math.AG/9805056.

 [SY] \name{B.\ Shelton} and \name{S.\  Yuzvinsky},  Koszul algebras from graphs and
hyperplane arrangements,
{\it J.\ London Math.\ Soc\/}.\ {\bf 56}  (1997), 477--490.

[Si] \name{D.\ Siersma}, Vanishing cycles and special fibres.\ 
Singularity theory and its applications.\ Part I (Coventry, 1988/1989), 
{\it Lecture Notes in Math\/}.\ {\bf 1462}  (1991),
292--301.

[ST] \name{D.\ Siersma} and \name{M.\  Tib\u ar}, Singularities at infinity and
their vanishing cycles II.\  Monodromy, {\it Publ.\ Res.\ Inst.\ Math.\
Sci\/}.\  {\bf 36}  (2000), 659--679.

[S] \name{E.\ H.\  Spanier}, {\it  Algebraic Topology\/}, McGraw-Hill
Book Co., New York, 1966.

[St] \name{J.\ Steenbrink}, Intersection form for quasi-homogeneous
singularities, {\it Compositio Math\/}.\ {\bf 34}  (1977), 211--223.

[V] \name{A.\ N.\  Varchenko}, Critical points of the product of powers of
linear functions and
families of bases of singular vectors, {\it Compositio Math\/}.\ {\bf
97}  (1995),
385--401.

[W] \name{G.\ W.\   Whitehead}, {\it Elements of Homotopy Theory\/}, 
{\it Grad.\ Texts in Math\/}.\ {\bf 61},\break Springer-Verlag, New York,
1978.

[Z] \name{O.\ Zariski}, On the problem of existence of algebraic
functions
of two variables possessing a given branch curve, {\it Amer.\ J.\
Math\/}.\ {\bf 51} 
(1929), 305--328.
\endreferences

\bye